# NEW MULTI-SAMPLE NONPARAMETRIC TESTS FOR PANEL COUNT DATA

By N. Balakrishnan and Xingqiu Zhao

*McMaster University and Hong Kong Polytechnic University*

This paper considers the problem of multi-sample nonparametric comparison of counting processes with panel count data, which arise naturally when recurrent events are considered. Such data frequently occur in medical follow-up studies and reliability experiments, for example. For the problem considered, we construct two new classes of nonparametric test statistics based on the accumulated weighted differences between the rates of increase of the estimated mean functions of the counting processes over observation times, wherein the nonparametric maximum likelihood approach is used to estimate the mean function instead of the nonparametric maximum pseudo-likelihood. The asymptotic distributions of the proposed statistics are derived and their finite-sample properties are examined through Monte Carlo simulations. The simulation results show that the proposed methods work quite well and are more powerful than the existing test procedures. Two real data sets are analyzed and presented as illustrative examples.

## 1. Introduction.

Consider a study that concerns some recurrent event, and suppose that each subject in the study gives rise to a counting process $N(t)$, denoting the total number of occurrences of the event of interest up to time $t$. Also suppose that for each subject, observations include only the values of $N(t)$ at discrete observation times or the numbers of occurrences of the event between the observation times. Such data are usually referred to as *panel count data* [Sun and Kalbfleisch (1995), Wellner and Zhang (2000)]. Our focus here will be on the situation when such a study involves $k$ $(\geq 2)$ groups. Let $\Lambda_l(t)$ denote the mean function of $N(t)$ corresponding to the $l$th group for $l = 1, \dots, k$. The problem of interest is then to test the hypothesis $H_0 : \Lambda_1(t) = \cdots = \Lambda_k(t)$.

A number of authors have discussed the analysis of recurrent event data when each subject in the study is observed continuously over an interval or when the exact times of occurrences of the recurrent event are known. For example, the book by Andersen et al. (1993) presents many of the commonly used statistical methods for the analysis of recurrent event data. In contrast, there exists limited research on the analysis of panel count data. Sun and Kalbfleisch (1995) and Wellner and Zhang (2000) studied estimation of the mean function of $N(t)$. Sun and Wei (2000) and Zhang (2002) discussed regression analysis for such data. To test the hypothesis $H_0$, Thall and Lachin (1988) suggested transforming the problem to a multivariate comparison problem and then applying a multivariate Wilcoxon-type rank test. Sun and Fang (2003) proposed a nonparametric procedure for this problem, but their procedure depends on the assumption that treatment indicators can be regarded as independent and identically distributed random variables, which may not be the case in practice. Park, Sun and Zhao (2007) proposed a class of nonparametric tests for the two-sample comparison based on the istonic regression estimator of the mean function of counting process. Zhang (2006) also presented nonparametric tests for the problem based on the nonparametric maximum pseudo-likelihood estimator, which is equivalent to the istonic regression estimator [Wellner and Zhang (2000)]. Also, Wellner and Zhang (2000) showed through Monte Carlo simulations that the nonparametric maximum likelihood estimator (NPMLE) of the mean function is more efficient than the nonparametric maximum pseudo-likelihood estimator (NPMPLE). However, no nonparametric tests have been discussed in the literature for panel count data based on the NPMLE, since the NPMLE is more complicated both theoretically and computationally. It is, therefore, particularly important to develop nonparametric tests based on the NPMLE for panel count data. One would naturally expect the tests based on the NPMLE to be more powerful than the tests based on the NPMPLE. However, unlike the isotonic regression estimate, the maximum likelihood estimate has no closed-form expression and its computation requires an iterative convex minorant algorithm. In this paper, we propose some nonparametric tests based on the maximum likelihood estimator and then compare them with the existing tests for the problem of multi-sample nonparametric comparison of counting processes with panel count data.

The rest of the paper is organized as follows. Section 2 discusses estimation of the mean function and the existing nonparametric tests for the hypothesis $H_0$ when only panel count data are available. The asymptotic normality of the functional of the NPMLE is established, while its proof is presented in Section 6. Section 3 presents two classes of nonparametric test statistics. The statistics, motivated by the property of the NPMLE and the idea used in survival analysis, are formulated as the integrated weighted difference between the rates of increase of the estimated mean functions corresponding to the pooled data and each group or two groups. The asymptotic normality of these test statistics is also established, while proofs are given in Section 6. In Section 4, finite-sample properties of the proposed test statistics are examined through Monte Carlo simulations. In Section 5, we apply the proposed methods to two data from a floating gallstones study and a bladder tumor study, respectively.

**2. Nonparametric maximum likelihood estimation of mean function.** Wellner and Zhang (2000) studied two estimators of the mean of a counting process with panel count data: the nonparametric maximum pseudo-likelihood estimator and the nonparametric maximum likelihood estimator. To describe the test statistics, we introduce first the NPMLE. Suppose that $N = \{N(t): t \geq 0\}$ is a nonhomogeneous Poisson process with the mean function $E(N(t)) = \Lambda_0(t)$. Suppose that $K$ is an integer-valued random variable and $T = \{T_{k,j}, j = 1, \ldots, k, k = 1, 2, \ldots\}$ is a random triangular array, where $T_{k,j-1} < T_{k,j}$ and $T_{k,0} = 0$, for $j = 1, \ldots, k$ and $k = 1, 2, \ldots$. We assume that $\{(K; T_{K,1}, \ldots, T_{K,K})\}$ are independent of $N$. Let $X = (K, T_K, N_K)$, where $T_k$ is the $k$th row of the triangular array $T$ and $N_k = (N(T_{k,1}), \ldots, N(T_{k,k}))$. Then, $X_i = (K_i, T_{K_i}, N_{i,K_i})$, $i = 1, \ldots, n$ is a random sample of size $n$ from the distribution of $X$. Let $\mathbf{X} = (X_1, \ldots, X_n)$. Then, the log-likelihood function for the mean function $\Lambda$ is

$$
\begin{aligned}
l_n(\Lambda|\mathbf{X}) = {} & \sum_{i=1}^{n}\sum_{j=1}^{K_i}(N_i(T_{K_i,j}) - N_i(T_{K_i,j-1}))\log(\Lambda(T_{K_i,j}) - \Lambda(T_{K_i,j-1})) \\
& - \sum_{i=1}^{n}\Lambda(T_{K_i,K_i})
\end{aligned}
$$

after omitting the parts independent of $\Lambda$.

Let $t_1 < \cdots < t_m$ denote the ordered distinct observation time points in the set of all observation time points $\{T_{K_i,j}, j = 1, \ldots, K_i, i = 1, \ldots, n\}$. Then the NPMLE of $\Lambda_0$, $\hat{\Lambda}_n$, is defined to be the nondecreasing, nonnegative step function with possible jumps only occurring at $t_\ell, \ell = 1, \ldots, m$ that maximizes $l_n(\Lambda|\mathbf{X})$. Wellner and Zhang (2000) gave the characteristic and the algorithm for computing this estimator, and studied its asymptotic properties.

Next, we need some more notation, some of which was introduced by Schick and Yu (2000) and Wellner and Zhang (2000). Let $\mathcal{B}$ denote the collection of Borel sets in $\mathcal{R}$, and let $\mathcal{B}_{[0,\tau]} = \{B \cap [0,\tau]: B \in \mathcal{B}\}$. Define measures $\mu_1$, $\mu_2$, $\mu_3$ and $\nu$ on $([0,\tau], \mathcal{B}_{[0,\tau]})$ by

$$
\mu_1(B) = \sum_{k=1}^{\infty} P(K = k) \sum_{j=1}^{k} P(T_{k,j} \in B | K = k),
$$

$$\mu_2(B_1 \times B_2) = \sum_{k=1}^{\infty} P(K=k) \sum_{j=1}^{k} P(T_{k,j-1} \in B_1, T_{k,j} \in B_2 | K=k),$$

$$\mu_3(B_1 \times B_2 \times B_3)$$
$$= \sum_{k=2}^{\infty} P(K=k) \sum_{j=1}^{k-1} P(T_{k,j-1} \in B_1, T_{k,j} \in B_2, T_{k,j+1} \in B_3 | K=k)$$

and

$$\nu(B_1 \times B_2) = \sum_{k=1}^{\infty} P(K=k) P(T_{k,k-1} \in B_1, T_{k,k} \in B_2 | K=k)$$

for $B, B_1, B_2, B_3 \in \mathcal{B}_{[0,\tau]}$.

The existing nonparametric tests [Park, Sun and Zhao (2007), Zhang (2006)] are based on the asymptotic normality of a smooth functional of the nonparametric maximum pseudo-likelihood estimator (the istonic regression estimator) $\tilde{\Lambda}_n$,

$$\int_0^{\tau} W(t)\{\tilde{\Lambda}_n(t) - \Lambda_0(t)\}\, d\mu_1(t) = P\left[\sum_{j=1}^{K} W(T_{K,j})\{\tilde{\Lambda}_n(T_{K,j}) - \Lambda_0(T_{K,j})\}\right],$$

where $W(t)$ is a weight function and $P$ is the probability measure of $X$, $Pf = \int f\, dP$. However, it is unknown if the asymptotic normality of the functional of the nonparametric maximum likelihood estimator $\int_0^{\tau} W(t)\{\hat{\Lambda}_n(t) - \Lambda_0(t)\}\, d\mu_1(t)$ still holds. We observe a key to the proof of such asymptotic normality is to use an important characteristic of the $\tilde{\Lambda}_n$ given by

$$\sum_{i=1}^{n} \sum_{j=1}^{K_i} \varphi(\tilde{\Lambda}_n(T_{K_i,j}))(\tilde{\Lambda}_n(T_{K_i,j}) - N_i(T_{K_i,j})) = 0 \tag{2.1}$$

for any real function $\varphi$. However, from (2.13) of Wellner and Zhang (2000), the corresponding characteristic of the NPMLE can be written as

$$\begin{aligned}\sum_{i=1}^{n} \Bigg[ \sum_{j=1}^{K_i-1} \hat{\Lambda}_n(T_{K_i,j}) \left\{ \frac{\Delta N_i(T_{K_i,j+1})}{\Delta \hat{\Lambda}_n(T_{K_i,j+1})} - \frac{\Delta N_i(T_{K_i,j})}{\Delta \hat{\Lambda}_n(T_{K_i,j})} \right\} \\ + \hat{\Lambda}_n(T_{K_i,K_i}) \left\{ 1 - \frac{\Delta N_i(T_{K_i,K_i})}{\Delta \hat{\Lambda}_n(T_{K_i,K_i})} \right\} \Bigg] = 0,\end{aligned} \tag{2.2}$$

where

$$\Delta\Lambda(T_{K,j}) = \Lambda(T_{K,j}) - \Lambda(T_{K,j-1})$$

and

$$\Delta N(T_{K,j}) = N(T_{K,j}) - N(T_{K,j-1}).$$

Equation (2.2) can be extended to the form

$$\sum_{i=1}^{n}\Bigg[\sum_{j=1}^{K_i-1}\varphi(\hat\Lambda_n(T_{K_i,j}))\hat\Lambda_n(T_{K_i,j})\bigg\{\frac{\Delta N_i(T_{K_i,j+1})}{\Delta\hat\Lambda_n(T_{K_i,j+1})}-\frac{\Delta N_i(T_{K_i,j})}{\Delta\hat\Lambda_n(T_{K_i,j})}\bigg\}$$

$$+\varphi(\hat\Lambda_n(T_{K_i,K_i}))\hat\Lambda_n(T_{K_i,K_i})\bigg\{1-\frac{\Delta N_i(T_{K_i,K_i})}{\Delta\hat\Lambda_n(T_{K_i,K_i})}\bigg\}\Bigg]=0, \tag{2.3}$$

which will be shown in Lemma 1. Clearly, the structure of (2.3) is different from that of (2.1) and is much more complicated. This is why the derivation of the asymptotic property of $\int_0^\tau W(t)\{\hat\Lambda_n(t)-\Lambda_0(t)\}\,d\mu_1(t)$ has not been done yet. So, we need to develop a new form of the test statistic when the NPMLE is used to estimate the mean function of counting process with panel count data. Motivated by such characteristic of the NPMLE, we define

$$f_\Lambda(X)=\sum_{j=1}^{K-1}W(T_{K,j})\Lambda(T_{K,j})\bigg\{\frac{\Delta\Lambda_0(T_{K,j+1})}{\Delta\Lambda(T_{K,j+1})}-\frac{\Delta\Lambda_0(T_{K,j})}{\Delta\Lambda(T_{K,j})}\bigg\}$$

$$+W(T_{K,K})\Lambda(T_{K,K})\bigg\{1-\frac{\Delta\Lambda_0(T_{K,K})}{\Delta\Lambda(T_{K,K})}\bigg\}. \tag{2.4}$$

It is easy to see that $Pf_\Lambda(X)$ can be expressed as

$$Pf_\Lambda(X)=\iiint W(u)\Lambda(u)\bigg\{\frac{\Lambda_0(v)-\Lambda_0(u)}{\Lambda(v)-\Lambda(u)}-\frac{\Lambda_0(u)-\Lambda_0(t)}{\Lambda(u)-\Lambda(t)}\bigg\}\,d\mu_3(t,u,v)$$

$$+\iint W(u)\Lambda(u)\bigg\{1-\frac{\Lambda_0(u)-\Lambda_0(t)}{\Lambda(u)-\Lambda(t)}\bigg\}\,d\nu(t,u).$$

For establishing asymptotic results on $Pf_{\hat\Lambda_n}(X)$, we need the following regularity conditions:

A. There exists a constant $K_0$ such that $P\{K\le K_0\}=1$ and that the random variables $T_{k,j}$'s take values in a bounded set $[\tau_0,\tau]$, where $\tau_0,\tau\in(0,\infty)$.

B. The mean function $\Lambda_0$ is strictly increasing such that $\Lambda_0(\tau_0)>0$ and $\Lambda_0(\tau)\le M$ for some constant $M\in(0,\infty)$.

C. There exists a constant $L_0$ such that

$$P\Big\{\min_{1\le j\le K}(\Lambda_0(T_{K,j})-\Lambda_0(T_{K,j-1}))\ge L_0\Big\}=1.$$

D. $E\{e^{cN(t)}\}$ is uniformly bounded for $t\in[0,\tau]$ and some constant $c$.

E. $\mu_1(\{\tau_0\})>0$ and for all $\tau_0<\tau_1<\tau_2<\tau$, $\mu_1((\tau_1,\tau_2))>0$.

Condition C holds if $\Lambda_0$ is differentiable, $\Lambda_0'$ has a positive lower bound in $[\tau_0,\tau]$ and $P\{\min_{1\le j\le K}(T_{K,j}-T_{K,j-1})\ge s_0\}=1$ for some fixed time $s_0$, where $s_0$ can be considered as the smallest length of consecutive observation times. Condition E holds if $P\{T_{K,1} = \tau_0\} > 0$ and $\mu_1'(t) > 0$ for $t \in (\tau_0, \tau)$.

Now, let $\Lambda_0^{-1}$ denote the inverse function of $\Lambda_0$, and let $W \circ \Lambda_0^{-1}$ denote composition of two functions $W$ and $\Lambda_0^{-1}$. Zhang (2006) established the asymptotic normality of $\int_0^\tau W(t)\{\tilde{\Lambda}_n(t) - \Lambda_0(t)\}\, d\mu_1(t)$, when $W \circ \Lambda_0^{-1}$ is not only bounded Lipschitz but also monotone. However, the assumption that $W \circ \Lambda_0^{-1}$ is monotone is not required for the tests with interval-censored data as a special case of panel count data [see Huang and Wellner (1995) and Zhang, Liu and Zhan (2001)]. Here, we do not need this monotonicity condition for $W \circ \Lambda_0^{-1}$.

THEOREM 2.1. *Suppose that Conditions* A, B, C, D *and* E *hold. Further, suppose that* $W(t)$ *is a bounded weight process such that* $W \circ \Lambda_0^{-1}$ *is a bounded Lipschitz function. Then, as* $n \to \infty$,

$$\sqrt{n} P f_{\hat{\Lambda}_n}(X) \longrightarrow U_w \tag{2.5}$$

*in distribution, where* $U_w$ *has a normal distribution with mean zero and variance* $\sigma_w^2$ *with*

$$\begin{aligned} \sigma_w^2 = E\Biggl[ & \sum_{j=1}^{K-1} W(T_{K,j})\Lambda_0(T_{K,j}) \biggl\{ \frac{\Delta N(T_{K,j+1})}{\Delta\Lambda_0(T_{K,j+1})} - \frac{\Delta N(T_{K,j})}{\Delta\Lambda_0(T_{K,j})} \biggr\} \\ & + W(T_{K,K})\Lambda_0(T_{K,K}) \biggl\{ 1 - \frac{\Delta N(T_{K,K})}{\Delta\Lambda_0(T_{K,K})} \biggr\} \Biggr]^2, \end{aligned} \tag{2.6}$$

*which can be consistently estimated by*

$$\begin{aligned} \hat{\sigma}_w^2 = \frac{1}{n} \sum_{i=1}^{n} \Biggl[ & \sum_{j=1}^{K_i-1} W(T_{K_i,j})\hat{\Lambda}_n(T_{K_i,j}) \biggl\{ \frac{\Delta N_i(T_{K_i,j+1})}{\Delta\hat{\Lambda}_n(T_{K_i,j+1})} - \frac{\Delta N_i(T_{K_i,j})}{\Delta\hat{\Lambda}_n(T_{K_i,j})} \biggr\} \\ & + W(T_{K_i,K_i})\hat{\Lambda}_n(T_{K_i,K_i}) \biggl\{ 1 - \frac{\Delta N_i(T_{K_i,K_i})}{\Delta\hat{\Lambda}_n(T_{K_i,K_i})} \biggr\} \Biggr]^2. \end{aligned} \tag{2.7}$$

## 3. Nonparametric tests.

Consider a longitudinal study that is concerned with some recurrent event and involves $n$ independent subjects, $n_l$ in the $l$th group with $n_1 + \cdots + n_k = n$ and $k \geq 2$. Let $N_i(t)$ denote the counting process arising from subject $i$ and $\Lambda_l(t)$ $(l = 1, \ldots, k)$ be as defined before, for $i = 1, \ldots, n$. Suppose that each subject is observed only at discrete time points $0 < T_{K_i,1} < \cdots < T_{K_i,K_i}$ and that no information is available about $N_i(t)$ between observation times; that is, only panel count data are available. For simplicity, assume that $H_0$ is true and let $\Lambda_0(t)$ denote the common mean function of the $N_i(t)$'s.

Let $\hat{\Lambda}_{n_l}$ denote the nonparametric maximum likelihood estimate of $\Lambda_l$ based on samples from all the subjects in the $l$th group, and $\hat{\Lambda}_n$ based on the pooled data. To test the hypothesis $H_0$, motivated by our asymptotic results in Section 2 and an idea commonly used in survival analysis [e.g., Andersen et al. (1993), Pepe and Fleming (1989), Petroni and Wolfe (1994), Cook, Lawless and Nadeau (1996), Zhang, Liu and Zhan (2001), Park, Sun and Zhao (2007), Zhang (2002, 2006)], we propose the statistics

$$
\begin{aligned}
U_n^{(l)} = \frac{1}{\sqrt{n}} \sum_{i=1}^{n} \Biggl[ \sum_{j=1}^{K_i-1} & W_n^{(l)}(T_{K_i,j}) \hat{\Lambda}_n(T_{K_i,j}) \\
& \times \left\{ \frac{\Delta\hat{\Lambda}_{n_l}(T_{K_i,j+1})}{\Delta\hat{\Lambda}_n(T_{K_i,j+1})} - \frac{\Delta\hat{\Lambda}_{n_l}(T_{K_i,j})}{\Delta\hat{\Lambda}_n(T_{K_i,j})} \right\} \\
& + W_n^{(l)}(T_{K_i,K_i}) \hat{\Lambda}_n(T_{K_i,K_i}) \left\{ 1 - \frac{\Delta\hat{\Lambda}_{n_l}(T_{K_i,K_i})}{\Delta\hat{\Lambda}_n(T_{K_i,K_i})} \right\} \Biggr]
\end{aligned} \tag{3.1}
$$

(for $l = 1, \ldots, k$) and

$$
\begin{aligned}
V_n^{(l)} = \frac{1}{\sqrt{n}} \sum_{i=1}^{n} \Biggl[ \sum_{j=1}^{K_i-1} & W_n^{(l)}(T_{K_i,j}) \hat{\Lambda}_n(T_{K_i,j}) \\
& \times \left\{ \left( \frac{\Delta\hat{\Lambda}_{n_1}(T_{K_i,j+1})}{\Delta\hat{\Lambda}_n(T_{K_i,j+1})} - \frac{\Delta\hat{\Lambda}_{n_1}(T_{K_i,j})}{\Delta\hat{\Lambda}_n(T_{K_i,j})} \right) \right. \\
& \quad \left. - \left( \frac{\Delta\hat{\Lambda}_{n_l}(T_{K_i,j+1})}{\Delta\hat{\Lambda}_n(T_{K_i,j+1})} - \frac{\Delta\hat{\Lambda}_{n_l}(T_{K_i,j})}{\Delta\hat{\Lambda}_n(T_{K_i,j})} \right) \right\} \\
& + W_n^{(l)}(T_{K_i,K_i}) \hat{\Lambda}_n(T_{K_i,K_i}) \\
& \times \left\{ \left( 1 - \frac{\Delta\hat{\Lambda}_{n_1}(T_{K_i,K_i})}{\Delta\hat{\Lambda}_n(T_{K_i,K_i})} \right) - \left( 1 - \frac{\Delta\hat{\Lambda}_{n_l}(T_{K_i,K_i})}{\Delta\hat{\Lambda}_n(T_{K_i,K_i})} \right) \right\} \Biggr]
\end{aligned} \tag{3.2}
$$

(for $l = 2, \ldots, k$), where $W_n^{(l)}(t)$'s are bounded weight processes. The statistic $U_n^{(l)}$ is the integrated weighted difference between the rates of increase of $\hat{\Lambda}_n$ and $\hat{\Lambda}_{n_l}$ over the observation times and the statistic $V_n^{(l)}$ has a similar meaning. For the selection of the weight process $W_n^{(l)}(t)$, a simple and natural choice is $W_n^{(1,l)}(t) = 1, l = 1, \ldots, k$. Another natural choice is $W_n^{(2,l)}(t) = Y_n(t) = \sum_{i=1}^{n} I(t \le T_{K_i,K_i})/n, l = 1, \ldots, k$, in which case weights are proportional to the number of subjects under observation. Based on groups, one may choose the weight process $W_n^{(l)}(t)$ as

$$
W_n^{(3,l)}(t) = Y_{n_l}(t) \quad \text{or} \quad \frac{Y_{n_l}(t)}{Y_n(t)} \quad \text{or} \quad \frac{Y_{n_1}(t) Y_{n_l}(t)}{Y_n(t)},
$$

where $Y_{n_l}(t)$ $(l = 1, \ldots, k)$ are defined as $Y_n(t)$, with the summation being only over subjects in the $l$th group. Some weight processes similar to $W_n^{(3)}$ have been used when recurrent event data are observed [see Andersen et al. (1993)]. In addition, $\sum_{i=1}^n I(t > T_{K_i,K_i})/n$ is also chosen as another weight process by Zhang (2006). Some other possible choices are

$$1 - Y_{n_l}(t), \qquad \frac{1 - Y_{n_l}(t)}{1 - Y_n(t)} \quad \text{and} \quad \frac{(1 - Y_{n_1}(t))(1 - Y_{n_l}(t))}{1 - Y_n(t)}.$$

Now, we state the asymptotic distribution of $\mathbf{U}_n = (U_n^{(1)}, \ldots, U_n^{(k)})^T$ and $\mathbf{V}_n = (V_n^{(2)}, \ldots, V_n^{(k)})^T$.

THEOREM 3.1. *Suppose that Conditions* A, B, C, D *and* E *hold. Further, suppose that* $W_n^{(l)}(t)$*'s are bounded weight processes and that there exists a bounded function* $W(t)$*, such that* $W \circ \Lambda_0^{-1}$ *is a bounded Lipschitz function and*

$$\left[\int_0^\tau \{W_n^{(l)}(t) - W(t)\}^2 \, d\mu_1(t)\right]^{1/2} = o_p(n^{-1/6}), \qquad l = 1, \ldots, k. \tag{3.3}$$

*Also, suppose that* $n_l/n \to p_l$ *as* $n \to \infty$*, where* $0 < p_l < 1$*,* $l = 1, \ldots, k$ *and* $p_1 + \cdots + p_k = 1$*. Then, under* $H_0 : \Lambda_1 = \cdots = \Lambda_k = \Lambda_0$*:*

(i) $\mathbf{U}_n$ *has an asymptotic normal distribution with mean vector* $\mathbf{0}$ *and covariance matrix*

$$\mathbf{\Sigma}_{\mathbf{U}_w} = \mathbf{\Gamma} \operatorname{diag}(\sigma_1^2, \sigma_2^2, \ldots, \sigma_k^2) \mathbf{\Gamma}', \tag{3.4}$$

*where*

$$\mathbf{\Gamma} = \begin{pmatrix} \sqrt{p_1} - \sqrt{\frac{1}{p_1}} & \sqrt{p_2} & \cdots & \sqrt{p_k} \\ \sqrt{p_1} & \sqrt{p_2} - \sqrt{\frac{1}{p_2}} & \cdots & \sqrt{p_k} \\ \cdots & \cdots & \cdots & \cdots \\ \sqrt{p_1} & \sqrt{p_2} & \cdots & \sqrt{p_k} - \sqrt{\frac{1}{p_k}} \end{pmatrix}$$

*and* $\sigma_1^2 = \cdots = \sigma_k^2 = \sigma_w^2$ *given in* (2.6).

(ii) $\mathbf{V}_n$ *has an asymptotic normal distribution with mean vector* $\mathbf{0}$ *and covariance matrix*

$$\mathbf{\Sigma}_{\mathbf{V}_w} = \mathbf{H} \operatorname{diag}(\sigma_1^2, \sigma_2^2, \ldots, \sigma_k^2) \mathbf{H}', \tag{3.5}$$

*where*

$$\mathbf{H} = \begin{pmatrix} -\sqrt{\frac{1}{p_1}} & \sqrt{\frac{1}{p_2}} & 0 & \cdots & 0 \\ -\sqrt{\frac{1}{p_1}} & 0 & \sqrt{\frac{1}{p_3}} & \cdots & 0 \\ \cdots & \cdots & \cdots & \cdots & \cdots \\ -\sqrt{\frac{1}{p_1}} & 0 & 0 & \cdots & \sqrt{\frac{1}{p_k}} \end{pmatrix}$$

*and* $\sigma_l^2$ *is as given in* (i).

(iii) *In addition, if*

$$\max_{1\le i\le n} E\left[\sum_{j=1}^{K_i}\{W_n^{(l)}(T_{K_i,j}) - W(T_{K_i,j})\}^2\right] \longrightarrow 0 \tag{3.6}$$

*for* $l = 1, \ldots, k$, *then* $\mathbf{\Sigma}_{\mathbf{U}_w}$ *and* $\mathbf{\Sigma}_{\mathbf{V}_w}$ *can be consistently estimated by*

$$\hat{\mathbf{\Sigma}}_{U_n} = \mathbf{\Gamma}_n \operatorname{diag}(\hat{\sigma}_1^2, \hat{\sigma}_2^2, \ldots, \hat{\sigma}_k^2)\mathbf{\Gamma}_n' \tag{3.7}$$

*and*

$$\hat{\mathbf{\Sigma}}_{V_n} = \mathbf{H}_n \operatorname{diag}(\hat{\sigma}_1^2, \hat{\sigma}_2^2, \ldots, \hat{\sigma}_k^2)\mathbf{H}_n', \tag{3.8}$$

*where*

$$\mathbf{\Gamma}_n = \begin{pmatrix} \sqrt{\frac{n_1}{n}} - \sqrt{\frac{n}{n_1}} & \sqrt{\frac{n_2}{n}} & \cdots & \sqrt{\frac{n_k}{n}} \\ \sqrt{\frac{n_1}{n}} & \sqrt{\frac{n_2}{n}} - \sqrt{\frac{n}{n_2}} & \cdots & \sqrt{\frac{n_k}{n}} \\ \cdots & \cdots & \cdots & \cdots \\ \sqrt{\frac{n_1}{n}} & \sqrt{\frac{n_2}{n}} & \cdots & \sqrt{\frac{n_k}{n}} - \sqrt{\frac{n}{n_k}} \end{pmatrix},$$

$$\mathbf{H}_n = \begin{pmatrix} -\sqrt{\frac{n}{n_1}} & \sqrt{\frac{n}{n_2}} & 0 & \cdots & 0 \\ -\sqrt{\frac{n}{n_1}} & 0 & \sqrt{\frac{n}{n_3}} & \cdots & 0 \\ \cdots & \cdots & \cdots & \cdots & \cdots \\ -\sqrt{\frac{n}{n_1}} & 0 & 0 & \cdots & \sqrt{\frac{n}{n_k}} \end{pmatrix}$$

*and*

$$\hat{\sigma}_l^2 = \frac{1}{n}\sum_{i=1}^{n}\left[\sum_{j=1}^{K_i-1} W_n^{(l)}(T_{K_i,j})\hat{\Lambda}_n(T_{K_i,j}) \right.$$

$$\times \left\{\frac{\Delta N_i(T_{K_i,j+1})}{\Delta\hat{\Lambda}_n(T_{K_i,j+1})} - \frac{\Delta N_i(T_{K_i,j})}{\Delta\hat{\Lambda}_n(T_{K_i,j})}\right\} \tag{3.9}$$

$$+ W_n^{(l)}(T_{K_i,K_i})\hat{\Lambda}_n(T_{K_i,K_i})\left\{1 - \frac{\Delta N_i(T_{K_i,K_i})}{\Delta\hat{\Lambda}_n(T_{K_i,K_i})}\right\}\Bigg]^2$$

*for* $l = 1, \ldots, k$.

Let $\mathbf{U}_0$ denote the first $(k-1)$ components of $\mathbf{U}_n$ and $\hat{\mathbf{\Sigma}}_0$ the matrix obtained by deleting the last row and column of $\hat{\mathbf{\Sigma}}_{U_n}$. Then, using Theorem 3.1, two tests can be carried out for testing $H_0$ by means of the statistic $\chi_0^2 = \mathbf{U}_0^T\hat{\mathbf{\Sigma}}_0^{-1}\mathbf{U}_0$ and $\mathbf{V}_n^T\hat{\mathbf{\Sigma}}_{V_n}^{-1}\mathbf{V}_n$, which have asymptotically a central $\chi^2$-distribution with $(k-1)$ degrees of freedom. This can be seen readily from the proof of the theorem.

Remark 1. If the weight process $W_n^{(l)}$ is symmetric about $X_1, \ldots, X_n$, then (3.6) is equivalent to

$$E\left[\sum_{j=1}^{K_1}\{W_n^{(l)}(T_{K_1,j}) - W(T_{K_1,j})\}^2\right] \longrightarrow 0.$$

Remark 2. For selection of weight processes, Zhang (2006) required that $W_n(t)$, $W(t)$ and $W \circ \Lambda_0^{-1}$ are monotone. These monotonicity assumptions restrict availability of weight processes. For example, the weight process $\frac{Y_{n_1}(t)Y_{n_2}(t)}{Y_{n_1}(t)+Y_{n_2}(t)}$ is often used in survival analysis, but it is not monotone. In addition, the monotonicity assumption on the weight process is not appropriate for deriving optimal tests under alternatives. In the above theorem, we have removed these assumptions. Therefore, compared to those stated in Zhang (2006), more weight processes are available here. It can be easily shown that the weight processes mentioned earlier all satisfy the conditions required by the theorem.

**4. Simulation study.** To examine the finite-sample properties of the proposed test statistics and compare them with those of the tests presented by Sun and Fang (2003), Park, Sun and Zhao (2007) and Zhang (2006), we carry out a simulation study for the two-sample comparison problem. When $k = 2$, the null hypothesis can be tested by $T_1 = U_n^{(1)}/\hat{\sigma}_U$ and $T_2 = V_n^{(2)}/\hat{\sigma}_V$, which have asymptotic standard normal distribution, where

$$\hat{\sigma}_U = \left\{\left(\sqrt{\frac{n_1}{n}} - \sqrt{\frac{n}{n_1}}\right)^2\hat{\sigma}_1^2 + \frac{n_2}{n}\hat{\sigma}_2^2\right\}^{1/2},$$

$$\hat{\sigma}_V = \left\{\frac{n}{n_1}\hat{\sigma}_1^2 + \frac{n}{n_2}\hat{\sigma}_2^2\right\}^{1/2},$$

and $U_n^{(1)}$, $V_n^{(2)}$ and $\hat{\sigma}_l$ are as given in (3.1), (3.2) and (3.9), respectively. Let $T_{\mathrm{SF}}$, $T_{\mathrm{PSZ}}$ and $T_{\mathrm{Z}}$ denote the tests presented by Sun and Fang (2003),

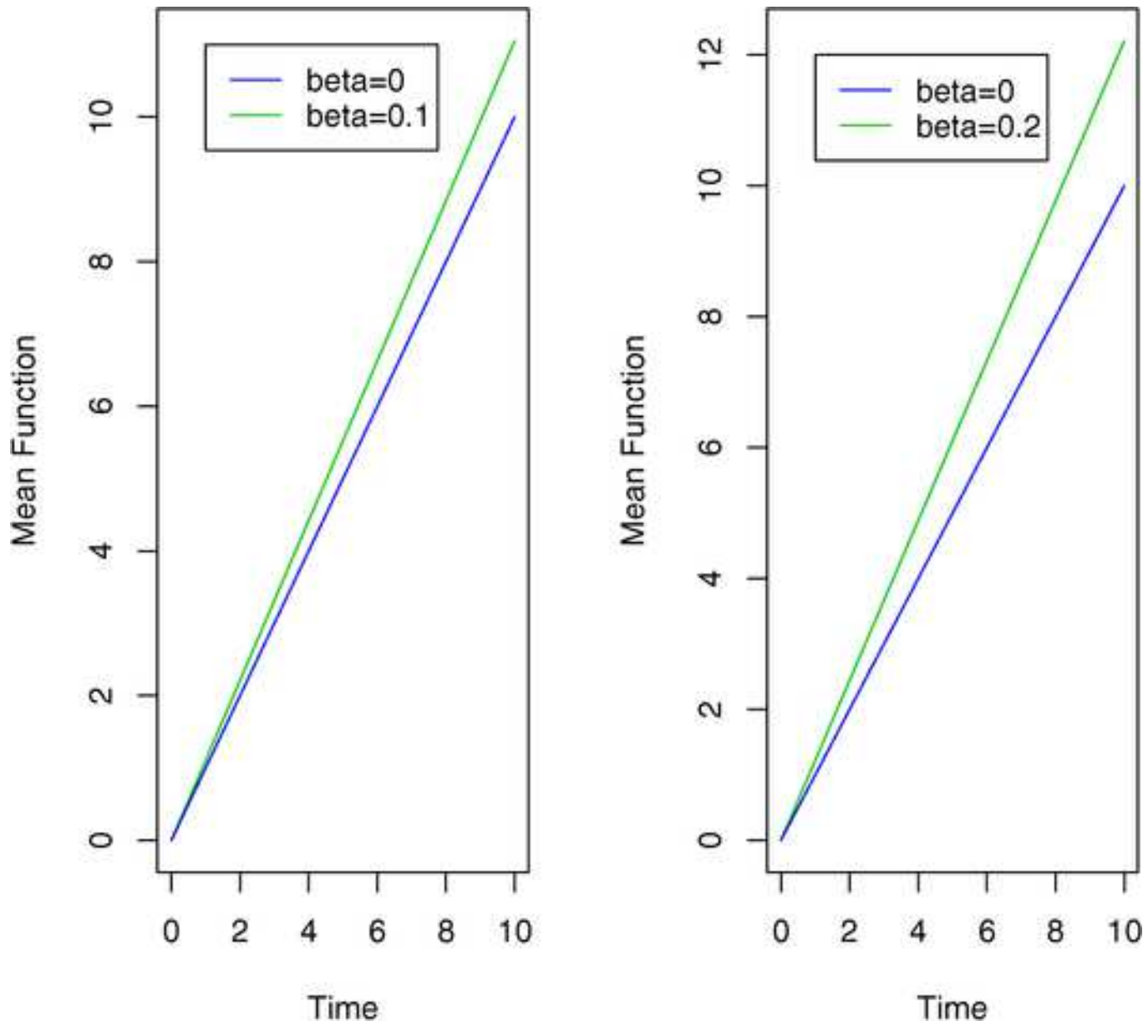

Fig. 1. *True mean functions for Case 1, with* $\nu = 1$ *and* $\beta = 0.1, 0.2$.

Park, Sun and Zhao (2007) and Zhang (2006), respectively. Here, we focus on evaluating the performance of $T_1$ and $T_2$ and comparing them to those of $T_{\mathrm{PSZ}}$, $T_{\mathrm{Z}}$ and $T_{\mathrm{SF}}$. Note that $T_{\mathrm{Z}} = T_{\mathrm{PSZ}}$ for $k = 2$. To generate panel count data $\{k_i, t_{ij}, n_{ij}, j = 1, \ldots, k_i, i = 1, \ldots, n\}$, we mimic medical follow-up studies such as the examples discussed in the next section. We first generate the number of observation times $k_i$ from the uniform distribution $U\{1, \ldots, 10\}$, and then, given $k_i$, we generate observation times $t_{ij}$'s from $U\{1, \ldots, 10\}$, for simplicity. To generate $n_{ij}$'s, we assume that $N_i$'s are nonhomogeneous Poisson or mixed Poisson processes. In particular, let $\{\nu_i, i = 1, \ldots, n\}$ be independent and identically distributed random variables, and given $\nu_i$, let $N_i(t)$ be a Poisson process with mean function $\Lambda_i(t|\nu_i) = E(N_i(t)|\nu_i)$. Let $S_l$ denote the set of indices for subjects in group $l$, $l = 1, 2$. For the objective of the study, we consider two cases as follows:

Case 1. $\Lambda_i(t|\nu_i) = \nu_i t$ for $i \in S_1$, $\Lambda_i(t) = \nu_i t \exp(\beta)$ for $i \in S_2$;

Case 2. $\Lambda_i(t|\nu_i) = \nu_i t$ for $i \in S_1$, $\Lambda_i(t) = \nu_i \sqrt{\beta t}$ for $i \in S_2$.

Figures 1 and 2 display the graphs of the true mean functions for two cases with $\nu = 1$ and different values of $\beta$. It can be seen that the two mean functions do not overlap in Case 1 and they cross over in Case 2.

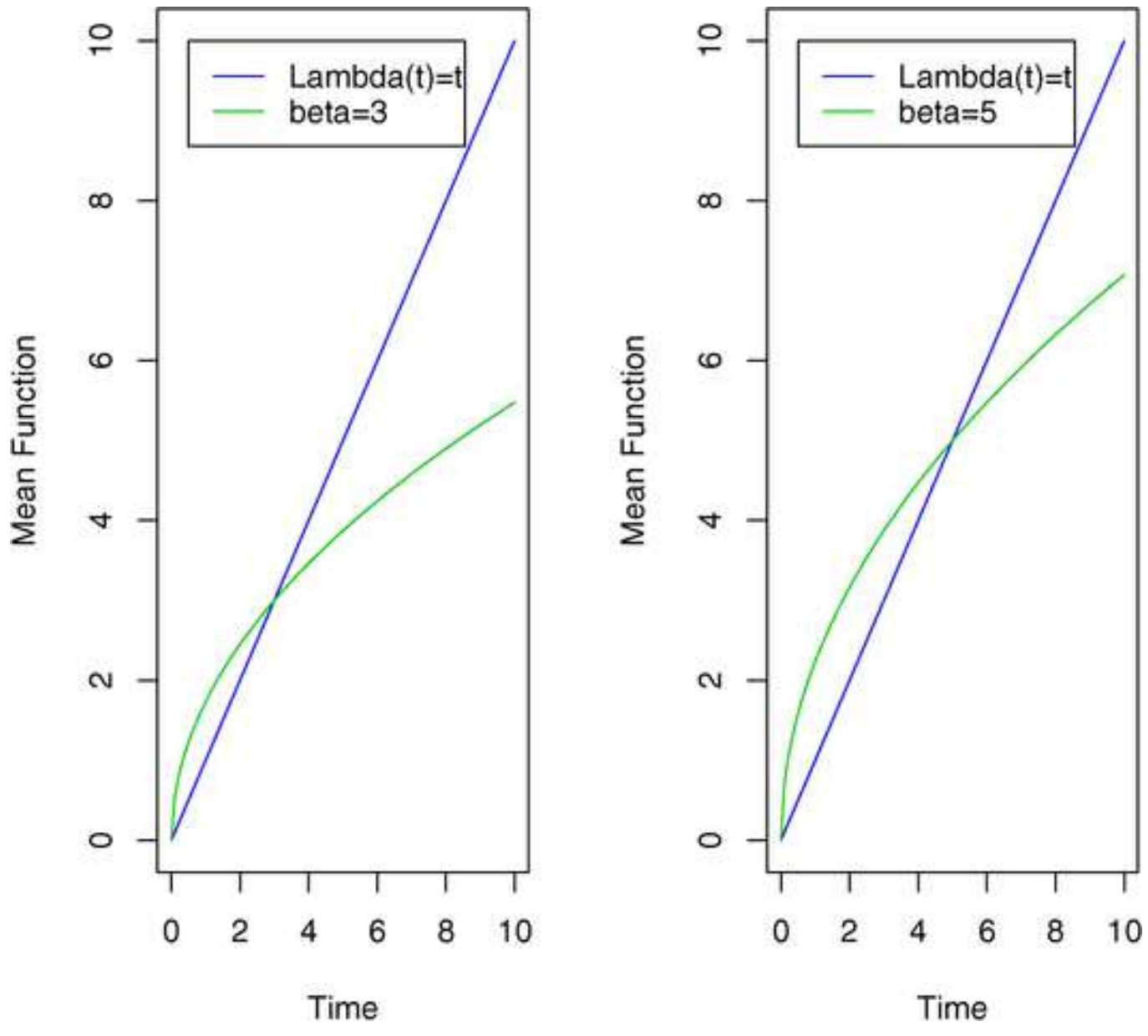

FIG. 2. *True mean functions for Case 2, with $\nu = 1$ and $\beta = 3, 5$.*

For each case, we consider $\nu_i = 1$ and $\nu_i \sim \text{Gamma}(2, 1/2)$, corresponding to Poisson and mixed Poisson processes, respectively. For each setting, we consider two sample sizes, $n_1 = n_2 = 50$ and 100, respectively. As mentioned earlier in Section 3, we choose the four weight processes

$$W_n^{(1)}(t) = 1, \qquad W_n^{(2)}(t) = Y_n(t) = \frac{1}{n}\sum_{i=1}^{n} I(t \le t_{k_i,k_i}),$$

$$W_n^{(3)}(t) = \frac{Y_{n_1}(t)Y_{n_2}(t)}{Y_n(t)} \quad \text{and} \quad W_n^{(4)}(t) = 1 - Y_n(t) = \frac{1}{n}\sum_{i=1}^{n} I(t > t_{k_i,k_i}).$$

The NPMLEs $\hat{\Lambda}_n$ and $\hat{\Lambda}_{n_l}$ are computed by using the modified iterative convex minorant algorithm (MICM) [see Wellner and Zhang (2000)]. All the results reported here are based on 1000 Monte Carlo replications using R software.

Tables 1–4 present the estimated sizes and powers of the proposed test statistics $T_1$ and $T_2$ and those of the test statistics $T_{\text{PSZ}}$, $T_{\text{Z}}$ and $T_{\text{SF}}$ [Park, Sun and Zhao (2007), Zhang (2006), Sun and Fang (2003)] at significance level $\alpha = 0.05$ for different values of $\beta$ and the four weight processes based on the simulated data for the two cases with $\nu_i = 1$ and $\nu_i \sim \text{Gamma}(2, 1/2)$, respectively. When $\nu_i = 1$, the $N_i(t)$'s are Poisson processes; when $\nu_i \sim \text{Gamma}(2, 1/2)$, the $N_i(t)$'s are mixed Poisson processes. The first part of the table is for the situation with the total sample size of 100, and the second part is for the situation with the total sample size of 200. For Case 1 considered here, the proposed tests display good power properties and the powers are close for the four weight processes. As expected, the power increases when the sample size increases, and the power decreases in the presence of more variability. As seen in Tables 1 and 2, the proposed tests with $W_n^{(1)}(t)$ have the best power performance, and the proposed tests based on the NPMLE are more powerful than the tests based on NPMPLE when more variability exists, as one would expect. For Case 2 considered here, the proposed tests also display good power properties, but the powers rely on choices of weight processes. As seen in Tables 3 and 4, the proposed tests with $W_n^{(4)}$ have the best power performance, and the proposed tests with appropriate weights based on NPMLE are much more powerful and more robust than those based on NPMPLE in this case. For example, when $\beta = 5, 8$ for mixed Poisson processes, the new tests with $W_n^{(4)}$ have good powers, but the tests $T_{\mathrm{PSZ}}$ and $T_{\mathrm{Z}}$ [Park, Sun and Zhao (2007), Zhang (2006)] with four weights and $T_{\mathrm{SF}}$ [Sun and Fang (2003)] have very poor powers. For all situations considered here, the performance of $T_1$ and $T_2$ are the same.

Note that the tests with different weights have different powers for Case 2. Let's explain why these results are reasonable. In this case, two true mean functions cross over at time $t = \beta$, the differences before this time point

Table 1

*Estimated size and power of the proposed test for Poisson processes in Case 1*

| | $T_2$ | | | | $T_{\mathrm{PSZ}}$ and $T_{\mathrm{Z}}$ | | | | |
|---|---|---|---|---|---|---|---|---|---|
| $\beta$ | $W_n^{(1)}$ | $W_n^{(2)}$ | $W_n^{(3)}$ | $W_n^{(4)}$ | $W_n^{(1)}$ | $W_n^{(2)}$ | $W_n^{(3)}$ | $W_n^{(4)}$ | $T_{\mathrm{SF}}$ |
| | | | | | $n_1 = n_2 = 50$ | | | | |
| 0.0 | 0.060 | 0.058 | 0.058 | 0.056 | 0.063 | 0.061 | 0.061 | 0.063 | 0.061 |
| 0.1 | 0.298 | 0.210 | 0.209 | 0.194 | 0.214 | 0.200 | 0.200 | 0.204 | 0.207 |
| 0.2 | 0.858 | 0.747 | 0.748 | 0.790 | 0.697 | 0.667 | 0.665 | 0.695 | 0.693 |
| 0.3 | 1.000 | 0.987 | 0.983 | 0.986 | 0.981 | 0.974 | 0.974 | 0.968 | 0.979 |
| | | | | | $n_1 = n_2 = 100$ | | | | |
| 0.0 | 0.047 | 0.047 | 0.047 | 0.049 | 0.044 | 0.046 | 0.046 | 0.045 | 0.043 |
| 0.1 | 0.542 | 0.472 | 0.471 | 0.489 | 0.423 | 0.405 | 0.405 | 0.411 | 0.422 |
| 0.2 | 0.993 | 0.967 | 0.964 | 0.991 | 0.958 | 0.948 | 0.947 | 0.952 | 0.950 |
| 0.3 | 1.000 | 1.000 | 1.000 | 1.000 | 1.000 | 1.000 | 1.000 | 1.000 | 1.000 |
| | $T_1$ | | | | $T_1$ | | | | |
| | $n_1 = n_2 = 50$ | | | | $n_1 = n_2 = 100$ | | | | |
| 0.0 | 0.052 | 0.051 | 0.051 | 0.053 | 0.051 | 0.049 | 0.049 | 0.050 | |
| 0.1 | 0.340 | 0.218 | 0.218 | 0.220 | 0.548 | 0.479 | 0.474 | 0.492 | |
| 0.2 | 0.868 | 0.787 | 0.764 | 0.798 | 0.996 | 0.976 | 0.974 | 0.993 | |
| 0.3 | 1.000 | 0.997 | 0.995 | 0.991 | 1.000 | 1.000 | 1.000 | 1.000 | |

TABLE 2
*Estimated size and power of the proposed test for mixed Poisson processes in Case 1*

| | $T_2$ | | | | $T_{PSZ}$ and $T_Z$ | | | | |
|---|---|---|---|---|---|---|---|---|---|
| $\beta$ | $W_n^{(1)}$ | $W_n^{(2)}$ | $W_n^{(3)}$ | $W_n^{(4)}$ | $W_n^{(1)}$ | $W_n^{(2)}$ | $W_n^{(3)}$ | $W_n^{(4)}$ | $T_{SF}$ |
| | | | | | $n_1 = n_2 = 50$ | | | | |
| 0.0 | 0.043 | 0.040 | 0.042 | 0.045 | 0.037 | 0.040 | 0.040 | 0.042 | 0.035 |
| 0.1 | 0.100 | 0.097 | 0.097 | 0.099 | 0.084 | 0.085 | 0.085 | 0.086 | 0.083 |
| 0.2 | 0.221 | 0.205 | 0.207 | 0.204 | 0.185 | 0.184 | 0.184 | 0.180 | 0.183 |
| 0.3 | 0.458 | 0.407 | 0.408 | 0.415 | 0.380 | 0.375 | 0.375 | 0.379 | 0.370 |
| | | | | | $n_1 = n_2 = 100$ | | | | |
| 0.0 | 0.043 | 0.041 | 0.041 | 0.046 | 0.048 | 0.045 | 0.045 | 0.044 | 0.046 |
| 0.1 | 0.140 | 0.125 | 0.125 | 0.138 | 0.114 | 0.111 | 0.111 | 0.112 | 0.111 |
| 0.2 | 0.410 | 0.364 | 0.362 | 0.368 | 0.317 | 0.307 | 0.307 | 0.314 | 0.316 |
| 0.3 | 0.708 | 0.663 | 0.662 | 0.672 | 0.596 | 0.592 | 0.592 | 0.593 | 0.590 |
| | | $T_1$ | | | | | $T_1$ | | |
| | | $n_1 = n_2 = 50$ | | | | | $n_1 = n_2 = 100$ | | |
| 0.0 | 0.054 | 0.048 | 0.048 | 0.046 | 0.048 | 0.047 | 0.047 | 0.051 | |
| 0.1 | 0.108 | 0.102 | 0.102 | 0.100 | 0.142 | 0.126 | 0.123 | 0.137 | |
| 0.2 | 0.216 | 0.205 | 0.206 | 0.207 | 0.412 | 0.388 | 0.390 | 0.391 | |
| 0.3 | 0.474 | 0.404 | 0.402 | 0.437 | 0.710 | 0.672 | 0.670 | 0.671 | |

TABLE 3
*Estimated power of the proposed test for Poisson processes in Case 2*

| | $T_2$ | | | | $T_{PSZ}$ and $T_Z$ | | | | |
|---|---|---|---|---|---|---|---|---|---|
| | $W_n^{(1)}$ | $W_n^{(2)}$ | $W_n^{(3)}$ | $W_n^{(4)}$ | $\beta W_n^{(1)}$ | $W_n^{(2)}$ | $W_n^{(3)}$ | $W_n^{(4)}$ | $T_{SF}$ |
| | | | | | $n_1 = n_2 = 50$ | | | | |
| 3 | 1.000 | 0.787 | 0.766 | 1.000 | 0.956 | 0.900 | 0.899 | 1.000 | 0.955 |
| 5 | 0.969 | 0.080 | 0.077 | 1.000 | 0.189 | 0.113 | 0.111 | 0.880 | 0.188 |
| 8 | 0.127 | 0.674 | 0.688 | 0.993 | 0.403 | 0.559 | 0.562 | 0.069 | 0.400 |
| | | | | | $n_1 = n_2 = 100$ | | | | |
| 3 | 1.000 | 0.964 | 0.960 | 1.000 | 0.999 | 0.993 | 0.993 | 1.000 | 0.999 |
| 5 | 1.000 | 0.078 | 0.079 | 1.000 | 0.290 | 0.140 | 0.139 | 0.988 | 0.284 |
| 8 | 0.222 | 0.935 | 0.939 | 1.000 | 0.670 | 0.843 | 0.846 | 0.082 | 0.667 |
| | | $T_1$ | | | | | $T_1$ | | |
| | | $n_1 = n_2 = 50$ | | | | | $n_1 = n_2 = 100$ | | |
| 3 | 1.000 | 0.784 | 0.769 | 1.000 | 1.000 | 0.958 | 0.955 | 1.000 | |
| 5 | 0.969 | 0.088 | 0.086 | 1.000 | 1.000 | 0.084 | 0.083 | 1.000 | |
| 8 | 0.130 | 0.675 | 0.689 | 0.995 | 0.232 | 0.932 | 0.935 | 1.000 | |

and after this time point have different signs; the difference after this point seems to dominate the difference before this point for the cases of $\beta = 3, 5$ and seems to be dominated by the difference before this point for the case of $\beta = 8$. When $\beta = 3, 5$, the tests with $W_n^{(1)}$ and $W_n^{(4)}$ have better powers than those with $W_n^{(2)}$ and $W_n^{(3)}$, and the test with $W_n^{(4)}$ has the largest power since it weights the difference at later times more than those with $W_n^{(1)}$, $W_n^{(2)}$ and $W_n^{(3)}$. In particular, when $\beta = 5$, the powers of the tests with $W_n^{(2)}$ and $W_n^{(3)}$ are very poor. This is because the small difference with large weights before this point and the large difference with small weights after this point seem to cancel each other. When $\beta = 8$, the tests with $W_n^{(2)}$, $W_n^{(3)}$ and $W_n^{(4)}$ perform better than the test with $W_n^{(1)}$. When $\beta = 8$, the biggest difference between two mean functions occurs at an earlier time, so that the tests with $W_n^{(2)}$ and $W_n^{(3)}$ have reasonable powers. But the test with $W_n^{(1)} = 1$ has a poor power though the difference at earlier times seems to dominate the difference at later times. This can be understood from the expressions of the test statistics $\mathbf{U}_n$ and $\mathbf{V}_n$, where the differences with different signs multiplied by the value of the mean function may cancel each other, since the mean function takes small values at earlier times and large values at later times. When $\beta = 8$, the test with $W_n^{(4)}$ still perform well. This is because it puts zero weight at earlier times and heavier weight at later times. Similar situations happened in real examples considered in the next section.

TABLE 4
*Estimated power of the proposed test for mixed Poisson processes in Case 2*

| | $\boldsymbol{T_2}$ | | | | $\boldsymbol{T_{\mathrm{PSZ}}}$ **and** $\boldsymbol{T_{\mathrm{Z}}}$ | | | | |
|---|---|---|---|---|---|---|---|---|---|
| | $\boldsymbol{W_n^{(1)}}$ | $\boldsymbol{W_n^{(2)}}$ | $\boldsymbol{W_n^{(3)}}$ | $\boldsymbol{W_n^{(4)}}$ | $\beta\boldsymbol{W_n^{(1)}}$ | $\boldsymbol{W_n^{(2)}}$ | $\boldsymbol{W_n^{(3)}}$ | $\boldsymbol{W_n^{(4)}}$ | $\boldsymbol{T_{\mathrm{SF}}}$ |
| | | | | | $n_1 = n_2 = 50$ | | | | |
| 3 | 0.858 | 0.301 | 0.294 | 0.992 | 0.386 | 0.321 | 0.318 | 0.708 | 0.380 |
| 5 | 0.424 | 0.078 | 0.078 | 0.943 | 0.089 | 0.071 | 0.071 | 0.289 | 0.086 |
| 8 | 0.062 | 0.255 | 0.263 | 0.771 | 0.117 | 0.158 | 0.158 | 0.039 | 0.111 |
| | | | | | $n_1 = n_2 = 100$ | | | | |
| 3 | 0.992 | 0.534 | 0.530 | 1.000 | 0.695 | 0.594 | 0.594 | 0.949 | 0.691 |
| 5 | 0.677 | 0.071 | 0.072 | 1.000 | 0.100 | 0.065 | 0.065 | 0.473 | 0.095 |
| 8 | 0.096 | 0.434 | 0.437 | 0.961 | 0.185 | 0.280 | 0.280 | 0.067 | 0.182 |
| | $T_1$ | | | | $T_1$ | | | | |
| | $n_1 = n_2 = 50$ | | | | $n_1 = n_2 = 100$ | | | | |
| 3 | 0.858 | 0.299 | 0.289 | 0.991 | 0.993 | 0.533 | 0.529 | 1.000 | |
| 5 | 0.396 | 0.074 | 0.074 | 0.942 | 0.685 | 0.068 | 0.066 | 1.000 | |
| 8 | 0.063 | 0.259 | 0.268 | 0.771 | 0.094 | 0.432 | 0.438 | 0.960 | |

To evaluate the asymptotic result given in Theorem 3.1, the quantile plots of the test statistic $T_2$ against the standard normal distribution are constructed. Figures 3 and 4 present the plots for the cases with $W_n(t) = W_n^{(1)}(t)$ and $n = 100$ and $n = 200$, respectively, and they clearly reveal that the asymptotic approximation is very good. Similar plots were obtained for test statistic $T_1$ and other situations as well.

In the above simulation study, we did examine all four weight processes suggested earlier in Section 3; in Case 1, the weight process $W_n^{(1)}$ yielded slightly higher power than the other three weight processes, and in Case 2, the weight process $W_n^{(4)}$ yielded the largest power. These simulation results suggest that, when the mean functions do not cross over, the test with the equal weight has a good power; otherwise, the test with the unequal and appropriate weight will also have a good power. In general, one can choose the weight process based on the behavior of the NPMLEs of the mean functions to improve power, since the true mean functions are unknown. When the difference of mean functions at earlier times dominate the difference at later times, the tests with $W_n^{(2)}$ and $W_n^{(3)}$ tend to have good powers; when the difference of mean functions at later times dominate the difference at earlier times, the test with $W_n^{(4)}$ tends to have a good power. In addition to the four processes considered here, some other weight processes can be found in Andersen et al. (1993), which discusses nonparametric treatment comparison based on recurrent event data. It would, therefore, be of great interest to investigate the problem of the selection of a weight process based on data.

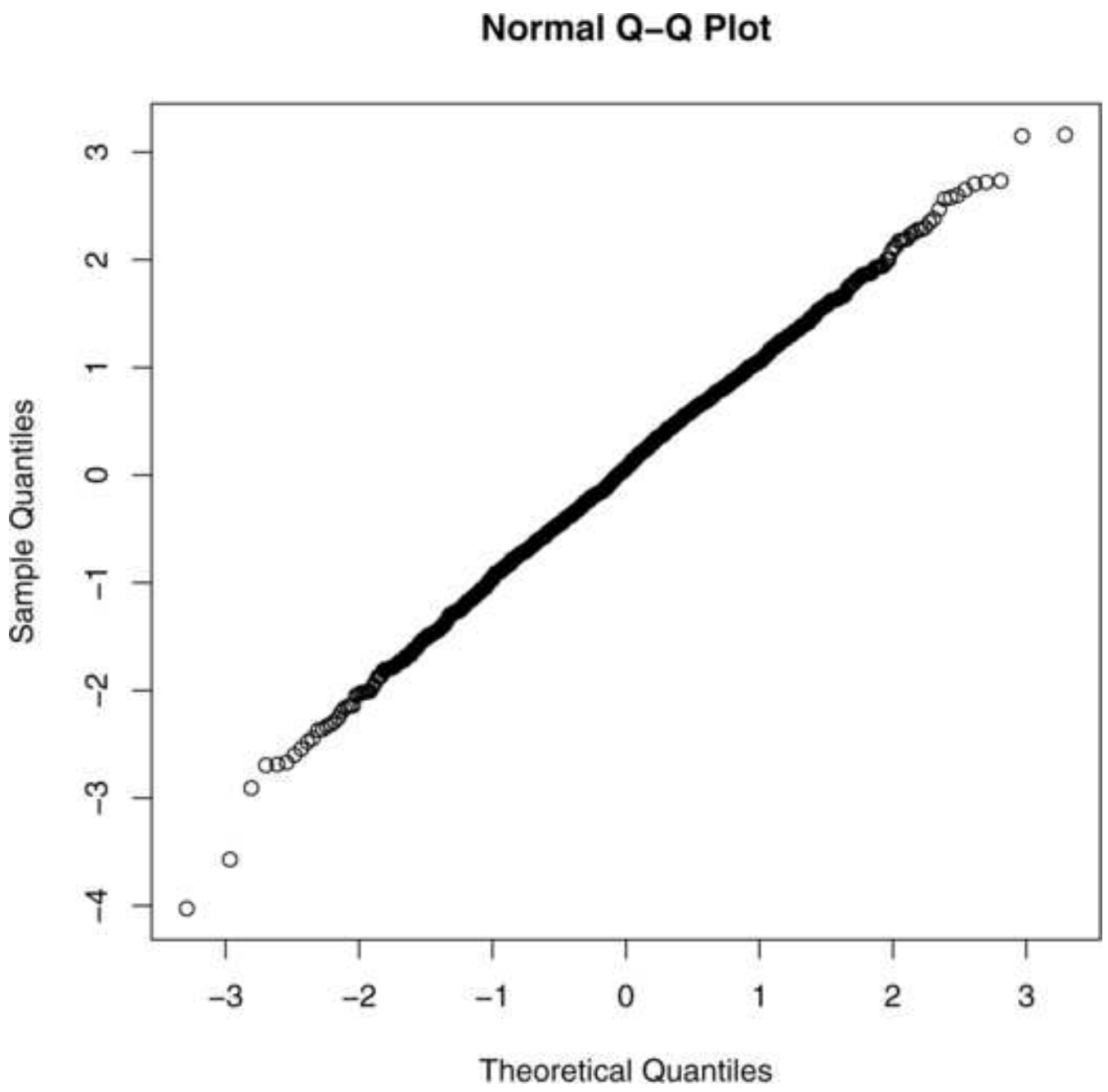

Fig. 3. *Simulation study. Normal quantile plot for $T_2$ ($n = 100$).*

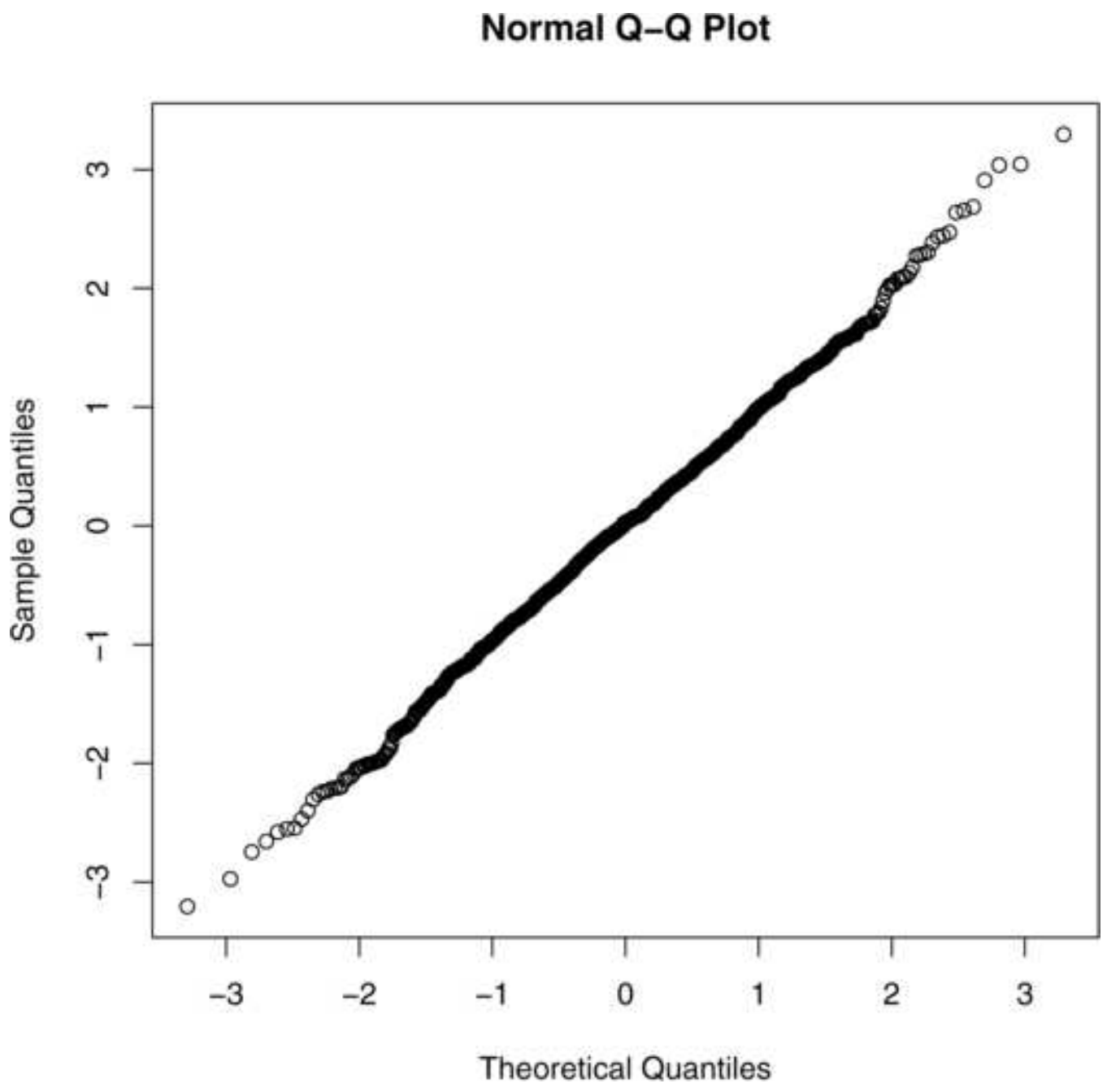

FIG. 4. *Simulation study. Normal quantile plot for $T_2$ $(n = 200)$.*

The new tests based on the NPMLE are more powerful and more robust than the existing tests based on the NPMPLE. One possible reason is that the NPMLE is more efficient than the NPMPLE. The main drawback of the NPMPLE is that the dependence of events within a subject is ignored. Another reason is that the structure of new test statistic is more reasonable, since it is based on the characteristic of the NPMLE.

**5. Illustrative examples.** To illustrate the proposed method, we consider here two examples: a floating gallstones study and a bladder tumor study.

5.1. *A floating gallstones study.* Thall and Lachin (1988) described a follow-up study on patients with floating gallstones. The data consist of the first year follow-up of the patients in two study groups, placebo (48) and high-dose chenodiol (65), from the National Cooperative Gallstone Study. The observed data include the successive visit times in study weeks and the associated counts of episodes of nausea for patients in different treatment groups [see Table 1 of Thall and Lachin (1988)]. The whole study consists of 916 patients who were randomized to placebo, low dose or high dose group and followed for up to two years. During the study, patients were scheduled to return for clinical visits at 1, 2, 3, 6, 9 and 12 months. In reality, most of the patients visited about six times within the first year, but actual visit times differ from patient to patient. Some patients had only one visit and some had 9 visits. As pointed out by Thall and Lachin (1988), there is no evidence that the number of observations and actual observation times are related to the incidence of nausea, and so it seems reasonable to assume that conditions required for the asymptotic results hold in this case. The problem of interest here is to compare the two treatment groups in terms of the incidence rates of nausea.

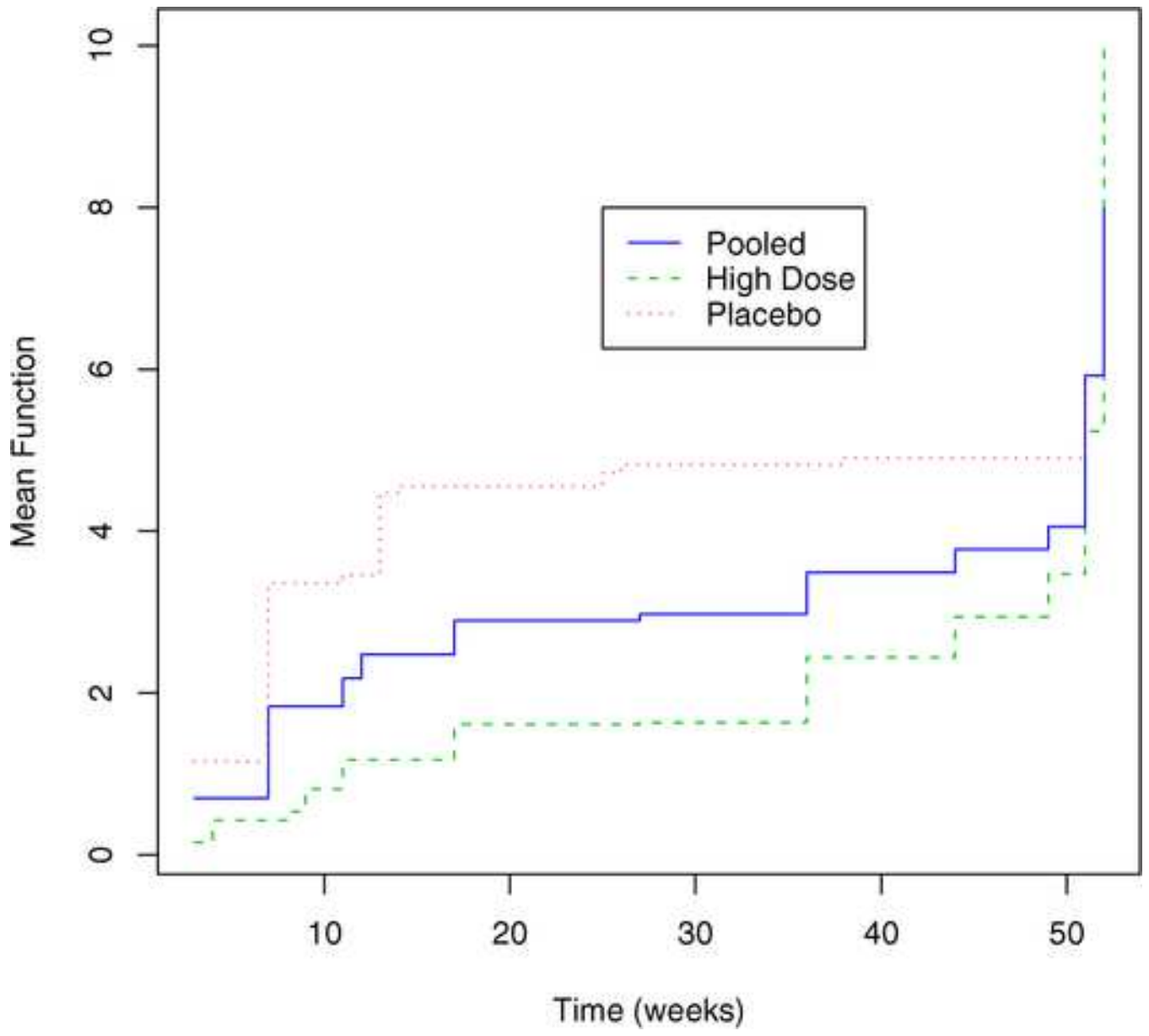

Fig. 5. *Floating gallstone study. Estimates of the mean functions.*

To test the difference between the two groups, we treated the placebo group as Group 1 ($\Lambda_1(t)$) and the high-dose chenodiol group as Group 2 ($\Lambda_2(t)$) and applied the proposed method to the data from 113 gallstone patients in the two groups to test the null hypothesis $H_0 : \Lambda_1(t) = \Lambda_2(t)$. The nonparametric maximum likelihood estimators of the incidence rates of nausea and the increments of the estimators are shown in Figures 5 and 6. We obtained $T_1 = 0.175$ and $T_2 = 0.206$ with $W_n(t) = W_n^{(1)}(t)$, giving $p$-values of 0.861 and 0.837 based on the standard normal distribution, and $T_1 = -337.221$, $-494.571$ and 241.159 and $T_2 = -193.238$, $-283.739$ and 138.311 with $W_n(t) = W_n^{(2)}(t)$, $W_n^{(3)}(t)$ and $1 - W_n^{(2)}(t)$, which correspond to $p$-values $\ll 0.0001$. The proposed tests with appropriate weights suggest that the incidence rates of nausea were significantly different for the patients in the two groups, and this agrees with the results given in Thall and Lachin (1988); the proposed unweighted test fails to reject $H_0$. This can be easily understood by looking at the behavior of increments of the estimators. From

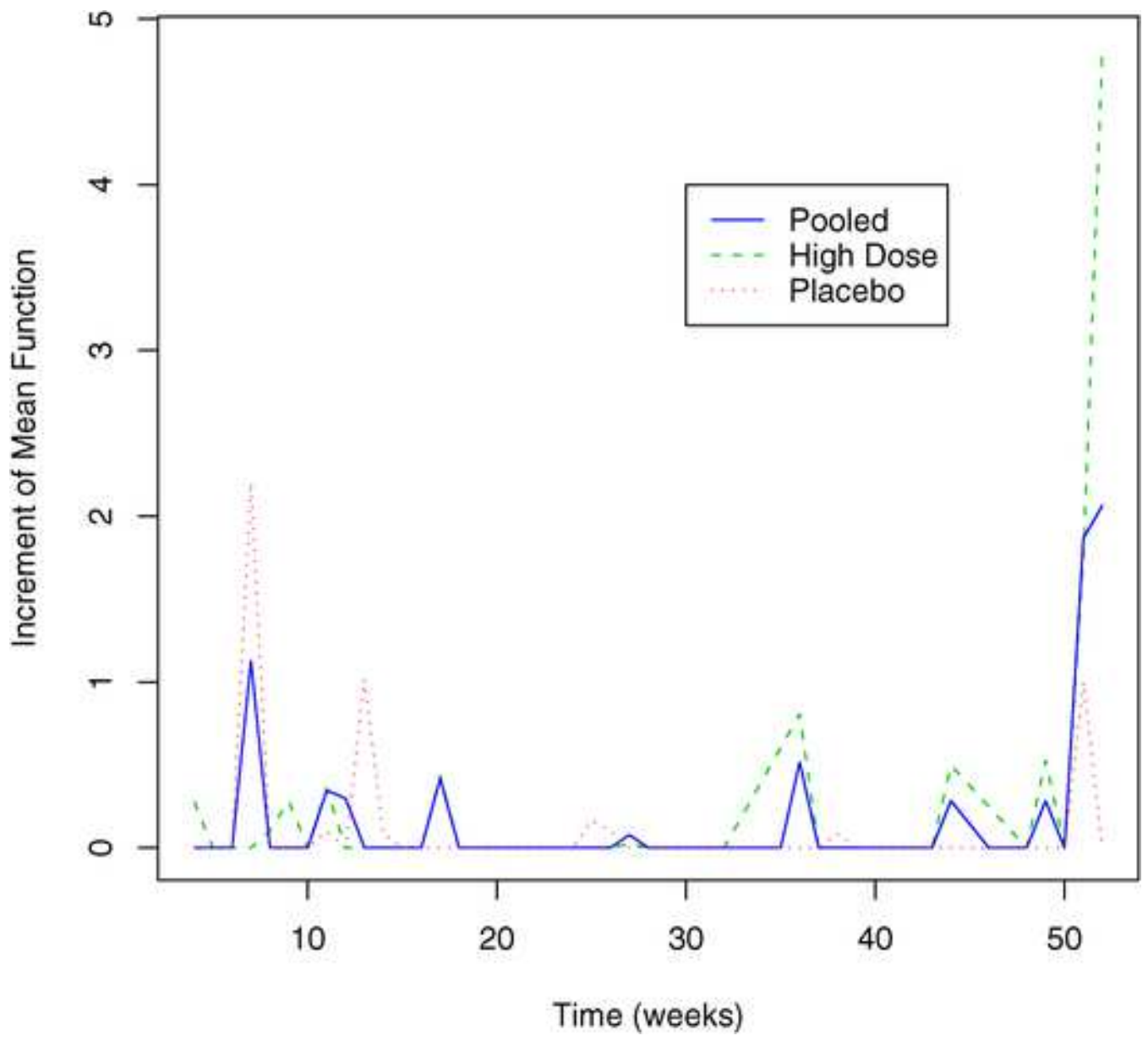

FIG. 6. *Floating gallstone study. Increments of the estimated mean functions.*

Figure 6, we can see clearly that the increment of the mean event rate in the placebo group is higher than that in the high dose group at earlier times and in contrast, the increment of the mean event rate in the high dose group is higher than that in the placebo group at later times in the year. So, the test with equal weights could not detect the difference between two groups. In comparison, the use of the approach in Sun and Fang (2003) gave a $p$-value of 0.1428; Park, Sun and Zhao (2007) gave $p$-values 0.454, 0.417 and 0.413 with three weights, respectively; and the tests presented by Zhang (2006) would give the same result as above. Thus, none of the existing tests based on NPMPLE can detect the difference of two treatments, and the proposed tests with suitable weights have detected successfully that, as we expected. One possible reason for this is that the nonparametric maximum likelihood estimator is more efficient than the nonparametric pseudo-likelihood estimator.

5.2. *A bladder tumor study.* We consider a bladder tumor study conducted by the Veterans Administration Co-operative Urological Research Group (VACURG), and the data are presented in Andrews and Herzberg (1985). For some earlier analyses of these data, one may refer to Byar, Blackard and The VACURG (1977), Byar (1980), Wellner and Zhang (2000), Sun and Wei (2000), and Zhang (2002, 2006). The data were obtained from a randomized clinical trial. All patients had superficial bladder tumors when they entered the trial, and they were assigned randomly to one of three treatments: placebo, thiotepa a pyridoxine. At subsequent follow-up visits, any tumors noticed were removed and treatment was continued. The study included 116 patients, of which there were 47 in placebo group, 38 in thiotepa group and 31 in pyridoxine. We can get a set of panel count data $\{k_i, t_{ij}, n_{ij}, j = 1, \ldots, k_i, i = 1, \ldots, n\}$ where for the $i$th patient, $k_i$ is the number of visits, $t_{ij}$'s are all visit times and $n_{ij}$ is total number of tumors until $t_{ij}$ $(j = 1, \ldots, k_i)$. The objective of the study is to determine the effect of treatment on the frequency of tumor recurrence.

Let $\Lambda_1(t), \Lambda_2(t)$ and $\Lambda_3(t)$ be the mean functions corresponding to the three treatment groups: placebo, thiotepa and pyridoxine, respectively. The nonparametric maximum likelihood estimators of mean functions and their increments from the three groups are presented in Figures 7 and 8, respectively. We observe from Figure 7 that the difference of the three groups becomes larger when the time increases. To test the null hypothesis $H_0 : \Lambda_1(t) = \Lambda_2(t) = \Lambda_3(t)$, we applied the proposed method to this panel count data. We obtained $\chi_0^2 = 3.617, 3.269$ and $p$-value $= 0.164, 0.195$ with $W_n(t) = 1$, $\chi_0^2 = 1196123, 300179.2$ and $p$-values $< 10^{-8}$ with $W_n(t) = Y_n(t)$, and $\chi_0^2 = 489000.4, 121908.1$ and $p$-values $< 10^{-8}$ with $W_n(t) = 1 - Y_n(t)$, based on the asymptotic distributions for test statistics $\mathbf{U}_n$ and $\mathbf{V}_n$ given in Theorem 3.1, respectively. The proposed tests having weights suggest that the frequency of tumor recurrence are significantly different for the patients in the three groups at 0.01 level of significance, while the proposed unweighted test fails to detect the difference. This can also be understood from the behavior of the increments of the estimated mean functions shown in Figure 8. Incidentally, through a regression analysis of the data from two treatments, placebo and thiotepa, Sun and Wei (2000) and Zhang (2002) concluded that thiotepa effectively reduces the recurrence of tumors. However, the existing test procedures [Sun and Fang (2003), Zhang (2006)] based on NPMPLE fail to reject the null hypothesis at level 0.05.

These examples illustrate that different weights may result in different conclusions, and the tests with appropriate weight process could lead to better power of the test. Therefore, the selection of a suitable weight process would be important for detecting difference between groups.

## 6. Proofs.

In this section we present the proofs of Theorems 2.1 and 3.1.

6.1. *Proof of Theorem 2.1.* We begin with some preliminary results. For convenience, let us first recall some notation given in Wellner and Zhang (2000). Set

$$\mathcal{F} = \{\Lambda : [0, \tau] \to [0, \infty) \mid \Lambda \text{ is nondecreasing}, \Lambda(0) = 0\}.$$

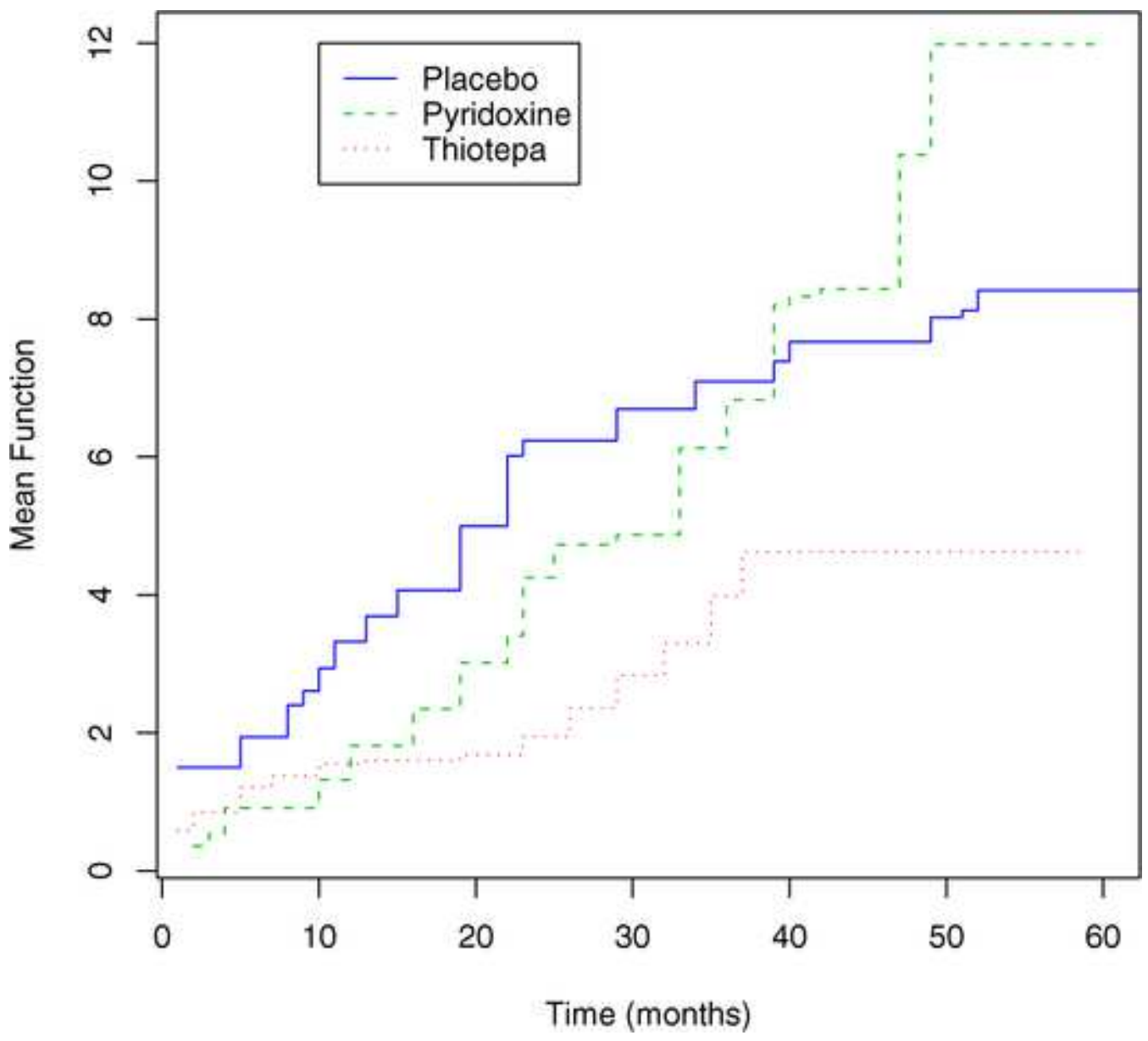

Fig. 7. *Bladder tumor study. Estimates of the mean functions.*

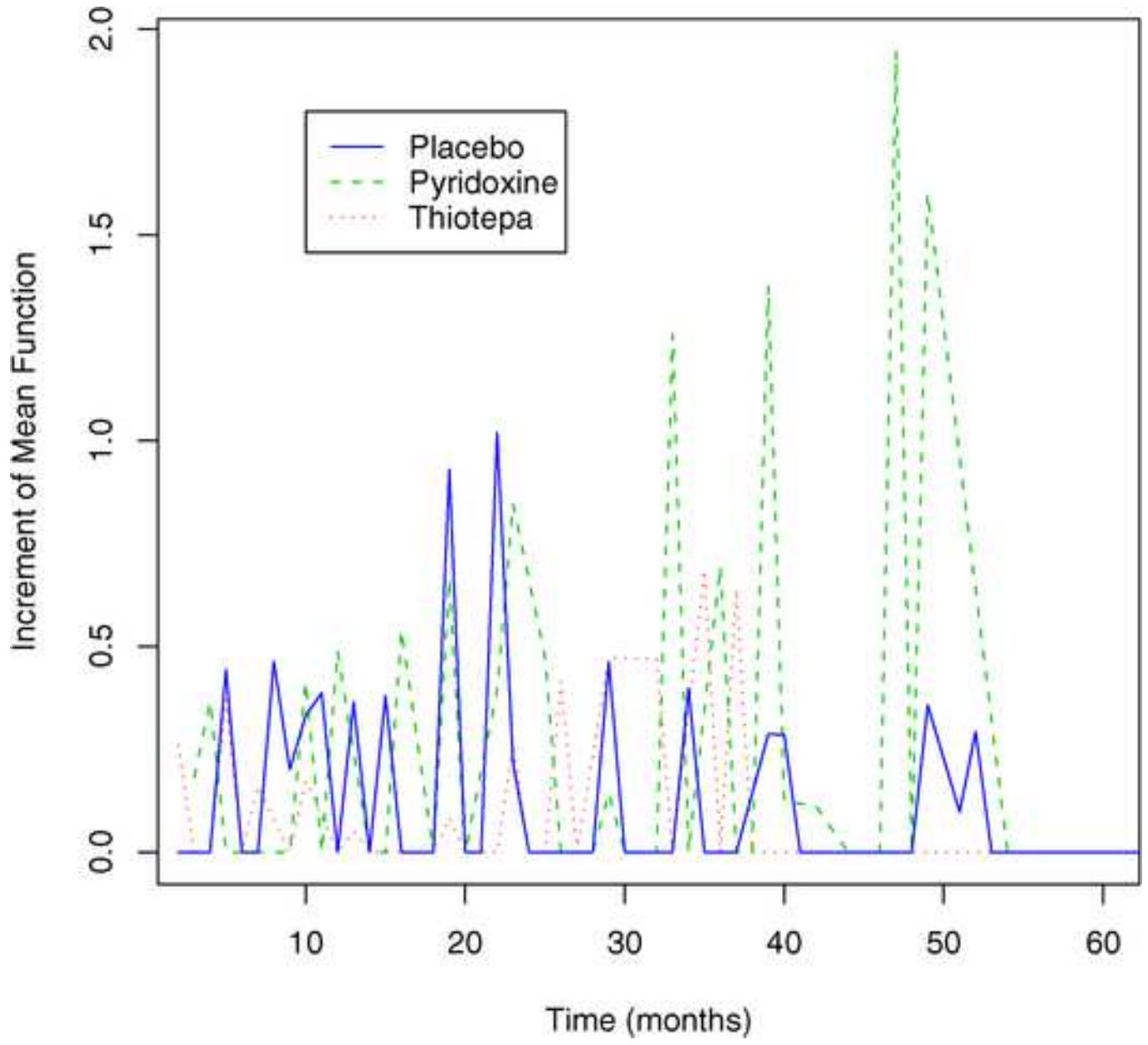

Fig. 8. *Bladder tumor study. Increments of the estimated mean functions.*

Let $t_1 < t_2 < \cdots < t_m$ denote the ordered distinct observation time points in the set of all observation time points $\{T_{K_i,j}, j = 1, \dots, K_i, i = 1, \dots, n\}$. Also, let $\Omega = \{\mathbf{u} = (u_1, u_2, \dots, u_m) : 0 \le u_1 \le \cdots \le u_m < \infty\}$ and the map $\mathcal{A} : \mathcal{F} \to \Omega$ be defined by

$$\mathbf{u} = \mathcal{A}(\Lambda) = (\Lambda(t_1), \Lambda(t_2), \dots, \Lambda(t_m)) \qquad \text{for all } \Lambda \in \mathcal{F}.$$

We also define a rank function $R$: $\{T_{K_i,j} : j = 1, 2, \dots, K_i; i = 1, 2, \dots, n\} \to \{1, 2, \dots, m\}$ such that

$$R(T_{K_i,j}) = s \qquad \text{if } T_{K_i,j} = t_s.$$

Then, the log-likelihood function can be rewritten as

$$\phi(\mathbf{u}|\mathbf{X}) = \sum_{i=1}^{n} \Bigg[ \sum_{j=1}^{K_i} \{N_i(T_{K_i,j}) - N_i(T_{K_i,j-1})\} \times \log\{u_{R(T_{K_i,j})} - u_{R(T_{K_i,j-1})}\} - u_{R(T_{K_i,K_i})} \Bigg],$$

and the NPMLE $\hat{\Lambda}_n$ of $\Lambda_0$ is then given by

$$(\hat{\Lambda}_n(t_1), \hat{\Lambda}(t_2), \dots, \hat{\Lambda}_n(t_m)) = \hat{\mathbf{u}}_n = \arg\max_{\mathbf{u} \in \Omega} \phi(\mathbf{u}|\mathbf{X}).$$

Set

$$\phi_\ell(\mathbf{u}) = \frac{\partial \phi(\mathbf{u}|\mathbf{X})}{\partial u_\ell} = \sum_{i=1}^{n} \phi_{i,\ell}(\mathbf{u}) \qquad \text{for } \ell = 1, 2, \dots, m,$$

where

$$\phi_{i,\ell}(\mathbf{u}) = \sum_{j=1}^{K_i - 1} \left\{ \frac{N_i(T_{K_i,j}) - N_i(T_{K_i,j-1})}{u_{R(T_{K_i,j})} - u_{R(T_{K_i,j-1})}} - \frac{N_i(T_{K_i,j+1}) - N_i(T_{K_i,j})}{u_{R(T_{K_i,j+1})} - u_{R(T_{K_i,j})}} \right\} 1_{\{T_{K_i,j} = t_\ell\}} + \left\{ \frac{N_i(T_{K_i,K_i}) - N_i(T_{K_i,K_i-1})}{u_{R(T_{K_i,K_i})} - u_{R(T_{K_i,K_i-1})}} - 1 \right\} 1_{\{T_{K_i,K_i} = t_\ell\}}.$$

LEMMA 1. *Let $\varphi$ be any real function. Then,*

$$\sum_{\ell=1}^{m} \varphi(\hat{u}_\ell) \left\{ \sum_{i=1}^{n} \phi_{i,\ell}(\hat{\mathbf{u}}) \right\} = 0. \tag{6.1}$$

PROOF. Let $\alpha_j = \hat{\Lambda}_n(t_j) - \hat{\Lambda}_n(t_{j-1})$, $j = 1, \ldots, m$. Using arguments similar to Proposition 2.1 of Groenebom (1996), we have

$$\sum_{j=i}^{m} \frac{\partial \phi(\hat{\mathbf{u}})}{\partial u_j} = 0 \qquad \text{if } \alpha_i > 0 \text{ or } i = 1.$$

Let $t_{k_1} < t_{k_2} < \cdots < t_{k_p}$ be jump points of $\hat{\Lambda}_n$. Then,

$$\sum_{\ell=k_j}^{m} \sum_{i=1}^{n} \phi_{i,\ell}(\hat{\mathbf{u}}) = 0, \qquad j = 1, \ldots, p,$$

and so

$$\sum_{k_j \le \ell < k_{j+1}} \sum_{i=1}^{n} \phi_{i,\ell}(\hat{\mathbf{u}}) = 0, \qquad j = 1, \ldots, p-1.$$

Thus,

$$\sum_{k_j \le \ell < k_{j+1}} \varphi(\hat{u}_\ell) \sum_{i=1}^{n} \phi_{i,\ell}(\hat{\mathbf{u}}) = 0, \qquad j = 1, \ldots, p-1,$$

since $\hat{u}_\ell = \hat{\Lambda}(t_\ell)$ is a constant for $k_j \le \ell < k_{j+1}$. Therefore, we conclude that

$$\sum_{\ell=1}^{m} \varphi(\hat{u}_\ell) \sum_{i=1}^{n} \phi_{i,\ell}(\hat{\mathbf{u}}) = 0.$$

Hence, the lemma follows. □

Now, let $\mu_i$ be as defined in Section 2, and let $d_i$ be the $L_2(\mu_i)$ metric on $\mathcal{F}$, $i = 1, 2$. Then, for $\Lambda_1, \Lambda_2 \in \mathcal{F}$,

$$\begin{aligned} d_1^2(\Lambda_1, \Lambda_2) &= \int |\Lambda_1(t) - \Lambda_2(t)|^2 \, d\mu_1(t) \\ &= E\left[\sum_{j=1}^{K} \{\Lambda_1(T_{K,j}) - \Lambda_2(T_{K,j})\}^2\right] \end{aligned} \tag{6.2}$$

and

$$\begin{aligned} &d_2^2(\Lambda_1, \Lambda_2) \\ &\quad = \iint |(\Lambda_1(s) - \Lambda_1(t)) - (\Lambda_2(s) - \Lambda_2(t))|^2 \, d\mu_2(s, t) \\ &\quad = E\left[\sum_{j=1}^{K} \{(\Lambda_1(T_{K,j}) - \Lambda_1(T_{K,j-1})) - (\Lambda_2(T_{K,j}) - \Lambda_2(T_{K,j-1}))\}^2\right]. \end{aligned} \tag{6.3}$$

If $P(K \le K_0) = 1$ for some constant $K_0$, then we have

$$\frac{1}{2} d_2(\Lambda_1, \Lambda_2) \le d_1(\Lambda_1, \Lambda_2) \le K_0 d_2(\Lambda_1, \Lambda_2). \tag{6.4}$$

Wellner and Zhang (2000) showed that

$$d_1(\hat{\Lambda}_n, \Lambda_0) \xrightarrow{\text{a.s.}} 0, \tag{6.5}$$

and hence that the uniform consistency of $\hat{\Lambda}_n$ can be shown by using arguments similar to Proposition 5 of Schick and Yu (2000) under Conditions A, B, D and E; that is,

$$\sup_{t \in [\tau_0, \tau]} |\hat{\Lambda}_n(t) - \Lambda_0(t)| \xrightarrow{\text{a.s.}} 0. \tag{6.6}$$

Note that the uniform consistency of $\hat{\Lambda}_n$ implies that for every $0 < \delta_0 < \min\{L_0/2, \Lambda_0(\tau_0)\}$ and any $\varepsilon > 0$, there exists a positive integer $N_\varepsilon$ such that

$$\sup_{n > N_\varepsilon} P\Big\{ \sup_{t \in [\tau_0, \tau]} |\hat{\Lambda}_n(t) - \Lambda_0(t)| > \delta_0 \Big\} < \varepsilon. \tag{6.7}$$

Here, we fix $\delta_0$. Let

$$\mathcal{F}_0 = \Big\{ \Lambda : \Lambda \in \mathcal{F}, \sup_{t \in [\tau_0, \tau]} |\Lambda(t) - \Lambda_0(t)| \le \delta_0 \Big\}. \tag{6.8}$$

Define $\hat{\Lambda}_n^*$ as

$$\hat{\Lambda}_n^* = \operatorname*{arg\,max}_{\Lambda \in \Omega \cap \mathcal{F}_0} \Bigg\{ \sum_{i=1}^{n} \sum_{j=1}^{K_i} (\Delta N_i(T_{K_i,j}) \log(\Delta\Lambda(T_{K_i,j})) - \Delta\Lambda(T_{K_i,j})) \Bigg\},$$

where $\Omega$ is the class of nondecreasing step functions with possible jumps only at the observation time points $\{T_{K_i,j}, j = 1, \dots, K_i, i = 1, \dots, n\}$. Clearly, we have

$$\sup_{n > N_\varepsilon} P(\hat{\Lambda}_n \ne \hat{\Lambda}_n^*) \le \sup_{n > N_\varepsilon} P\Big\{ \sup_{t \in [\tau_0, \tau]} |\hat{\Lambda}_n(t) - \Lambda_0(t)| > \delta_0 \Big\} < \varepsilon. \tag{6.9}$$

LEMMA 2. *We have* $d_1(\hat{\Lambda}_n^*, \Lambda_0) = O_p(n^{-1/3})$.

PROOF. To establish the rate of convergence for $\hat{\Lambda}_n^*$, we shall apply Theorem 3.2.5 of Van der Vaart and Wellner (1996). Define

$$\begin{aligned} m_\Lambda(X) = \sum_{j=1}^{K} [&(N(T_{K,j}) - N(T_{K,j-1})) \log\{\Lambda(T_{K,j}) - \Lambda(T_{K,j-1})\} \\ &- \{\Lambda(T_{K,j}) - \Lambda(T_{K,j-1})\}] \end{aligned} \tag{6.10}$$

and

$$\mathbb{M}(\Lambda) = Pm_\Lambda(X). \tag{6.11}$$

Let $h(x) = x(\log(x) - 1) + 1$. Then, $h(x) \geq \frac{1}{5}(x-1)^2$ for $x$ in a neighborhood of $x = 1$. Thus, in a neighborhood of $\Lambda_0$,

$$\begin{aligned}
&\mathbb{M}(\Lambda_0) - \mathbb{M}(\Lambda) \\
&\quad = \mathrm{P}\left[\sum_{j=1}^{K}\{\Lambda(T_{K,j}) - \Lambda(T_{K,j-1})\} h\left(\frac{\Lambda_0(T_{K,j}) - \Lambda_0(T_{K,j-1})}{\Lambda(T_{K,j}) - \Lambda(T_{K,j-1})}\right)\right] \\
&\quad = \iint \{\Lambda(u) - \Lambda(v)\} h\left(\frac{\Lambda_0(u) - \Lambda_0(v)}{\Lambda(u) - \Lambda(v)}\right) d\mu_2(u,v) \\
&\quad \geq \frac{1}{5}\iint \{\Lambda(u) - \Lambda(v)\}\left\{\frac{\Lambda_0(u) - \Lambda_0(v)}{\Lambda(u) - \Lambda(v)} - 1\right\}^2 d\mu_2(u,v) \\
&\quad = \frac{1}{5}\iint \frac{\{(\Lambda_0(u) - \Lambda(u)) - (\Lambda_0(v) - \Lambda(v))\}^2}{\Lambda(u) - \Lambda(v)} d\mu_2(u,v) \\
&\quad \geq c_1 d_1^2(\Lambda, \Lambda_0)
\end{aligned}$$

for some constant $c_1$, and hence the separation condition of the theorem is satisfied. Also, let

$$\mathcal{F}_\delta = \{\Lambda : d_1(\Lambda, \Lambda_0) \leq \delta, \Lambda \in \mathcal{F}_0\} \qquad (\delta > 0) \tag{6.12}$$

and

$$\mathcal{M}_\delta = \{m_\Lambda(X) - m_{\Lambda_0}(X) : \Lambda \in \mathcal{F}_\delta\}. \tag{6.13}$$

Note that $\mathcal{F}_\delta$ is a class of monotone nondecreasing functions. Then, it follows from Theorem 2.7.5 of Van der Vaart and Wellner (1996) that for any $\eta > 0$, there exists a set of brackets $\{[\Lambda_i^L, \Lambda_i^R] : i = 1, \ldots, J\}$, where $J \leq e^{c_2/\eta}$ for some constant $c_2$ and $d_1(\Lambda_i^L, \Lambda_i^R) \leq \eta$ such that for any $\Lambda \in \mathcal{F}_\delta$, $\Lambda_i^L \leq \Lambda \leq \Lambda_i^R$ for some $i$ with $1 \leq i \leq J$. Note that $\Lambda_i^L, \Lambda_i^R$ $(i = 1, \ldots, J)$ may not belong to $\mathcal{F}_\delta$, and so they may not have a uniform positive lower bound and a uniform finite upper bound in $[\tau_0, \tau]$. Also note that for any $\Lambda \in \mathcal{F}_0$, we have from Conditions A, B and C that

$$0 < \Lambda_0(\tau_0) - \delta_0 \leq \Lambda_0(t) - \delta_0 \leq \Lambda(t) \leq \Lambda_0(t) + \delta_0 \leq \Lambda_0(\tau) + \delta_0 \leq M + \delta_0$$

for $t \in [\tau_0, \tau]$ and

$$\begin{aligned}
0 < L_0 - 2\delta_0 &\leq \Delta\Lambda_0(T_{K,j}) - 2\delta_0 \\
&\leq \Delta\Lambda(T_{K,j}) \leq \Delta\Lambda_0(T_{K,j}) + 2\delta_0 \leq 2M + 2\delta_0
\end{aligned}$$

for $j=1,\dots,K$ with probability 1. Hence, for $\mathcal{M}_\delta$, we can construct a set of brackets $\{[M_i^L(X),M_i^R(X)]: i=1,\dots,J\}$ as follows:

$$\begin{aligned}M_i^L(X)=\sum_{j=1}^K&[\Delta N(T_{K,j})\\&\times\log\{\max(\Lambda_i^L(T_{K,j})-\Lambda_i^R(T_{K,j-1}),\Delta\Lambda_0(T_{K,j})-2\delta_0)\}\\&\qquad-\{\Lambda_i^R(T_{K,j})-\Lambda_i^L(T_{K,j-1})\}]\\&-m_{\Lambda_0}(X)\end{aligned}$$

and

$$\begin{aligned}M_i^R(X)=\sum_{j=1}^K&[\Delta N(T_{K,j})\log\{\Lambda_i^R(T_{K,j})-\Lambda_i^L(T_{K,j-1})\}\\&-\{\Lambda_i^L(T_{K,j})-\Lambda_i^R(T_{K,j-1})\}]-m_{\Lambda_0}(X).\end{aligned}$$

Set $\|\cdot\|_{P,B}$ be the Bernstein norm as defined in Van der Vaart and Wellner (1996) and $N_{[\cdot]}$ the braking number for the class $M_\delta$. Then, it follows from Condition D that

$$\|M_i^R(X)-M_i^L(X)\|_{P,B}^2\le c_3d_1^2(\Lambda_i^L,\Lambda_i^R)\le c_3\eta^2,\qquad i=1,\dots,J$$

for some constant $c_3$ and for any $\Lambda\in\mathcal{F}_\delta$, and

$$\|m_\Lambda(X)-m_{\Lambda_0}(X)\|_{P,B}^2\le c_4d_1^2(\Lambda,\Lambda_0)\le c_4\delta^2$$

for some constant $c_4$. So,

$$\log N_{[\cdot]}(\eta,\mathcal{M}_\delta,\|\cdot\|_{P,B})\le c_5\eta^{-1}$$

for some constant $c_5$. Hence, by applying Lemma 3.4.3 of Van der Vaart and Wellner (1996), we have

$$E^*\|\sqrt{n}(P_n-P)\|_{\mathcal{M}_\delta}\le c_6\phi_n(\delta)$$

for some constant $c_6$, where $E^*$ denotes the outer expectation, $P_n$ is the empirical measure corresponding to $X$, $P_nf=\sum_{i=1}^n f(X_i)/n$ and $\phi_n(\delta)=\delta^{1/2}+\delta^{-1}n^{-1/2}$. Now, upon using Theorem 3.2.5 of Van der Vaart and Wellner (1996), $d_1(\hat\Lambda_n^*,\Lambda_0)$ converges in probability to zero of order at least $n^{-1/3}$. This completes the proof of the lemma. □

Now we turn to the proof of Theorem 2.1. First, note that

$$\sqrt{n}\mathrm{P}f_{\hat\Lambda_n}(X)=-I_{1n}+I_{2n}+I_{3n},\tag{6.14}$$

where

$$I_{1n} = \sqrt{n}(P_n - P)f_{\hat{\Lambda}_n}(X),$$

$$I_{2n} = \sqrt{n}P_n\Bigg[\sum_{j=1}^{K-1} W(T_{K,j})\hat{\Lambda}_n(T_{K,j})\bigg\{\frac{\Delta N(T_{K,j+1})}{\Delta\hat{\Lambda}_n(T_{K,j+1})} - \frac{\Delta N(T_{K,j})}{\Delta\hat{\Lambda}_n(T_{K,j})}\bigg\} + W(T_{K,K})\hat{\Lambda}_n(T_{K,K})\bigg\{1 - \frac{\Delta N(T_{K,K})}{\Delta\hat{\Lambda}_n(T_{K,K})}\bigg\}\Bigg]$$

and

$$I_{3n} = \sqrt{n}P_n\Bigg[\sum_{j=1}^{K-1} W(T_{K,j})\hat{\Lambda}_n(T_{K,j})\bigg\{\frac{\Delta\Lambda_0(T_{K,j+1}) - \Delta N(T_{K,j+1})}{\Delta\hat{\Lambda}_n(T_{K,j+1})} - \frac{\Delta\Lambda_0(T_{K,j}) - \Delta N(T_{K,j})}{\Delta\hat{\Lambda}_n(T_{K,j})}\bigg\} + W(T_{K,K})\hat{\Lambda}_n(T_{K,K})\frac{\Delta N(T_{K,K}) - \Delta\Lambda_0(T_{K,K})}{\Delta\hat{\Lambda}_n(T_{K,K})}\Bigg].$$

Let

$$g_\Lambda(X) = \sum_{j=1}^{K-1} W(T_{K,j})\Lambda(T_{K,j})\bigg\{\frac{\Delta\Lambda_0(T_{K,j+1}) - \Delta N(T_{K,j+1})}{\Delta\Lambda(T_{K,j+1})} - \frac{\Delta\Lambda_0(T_{K,j}) - \Delta N(T_{K,j})}{\Delta\Lambda(T_{K,j})}\bigg\} + W(T_{K,K})\Lambda(T_{K,K})\frac{\Delta N(T_{K,K}) - \Delta\Lambda_0(T_{K,K})}{\Delta\Lambda(T_{K,K})}.$$

Note that

$$I_{3n} = \sqrt{n}(P_n - P)g_{\hat{\Lambda}_n}(X) = I_{4n} + I_{5n},$$

where

$$I_{4n} = \sqrt{n}(P_n - P)\{g_{\hat{\Lambda}_n}(X) - g_{\Lambda_0}(X)\}$$

and

$$I_{5n} = \sqrt{n}(P_n - P)g_{\Lambda_0}(X).$$

It is easy to see that $I_{5n}$ is a $U$-statistic and has an asymptotic normal distribution with mean zero and variance $\sigma_w^2$ that can be consistently estimated by $\hat{\sigma}_w^2$ as given in the statement of the theorem. Hence, it is sufficient to show that $I_{1n}$, $I_{2n}$ and $I_{4n}$ all converge in probability to zero.

We will show the convergence of $I_{1n}$ first. Let $I_{1n}^*$ denote the version of $I_{1n}$ obtained by replacing $\hat{\Lambda}_n$ with $\hat{\Lambda}_n^*$. Then, to prove that $I_{1n}$ converges to zero in probability, it is sufficient to show that $I_{1n}^* = o_p(1)$, since $P\{\hat{\Lambda}_n \neq \hat{\Lambda}_n^*\} < \varepsilon$. Let

$$\mathcal{F}_1 = \{f_\Lambda(X) : \Lambda \in \mathcal{F}_0\}.$$

Also, let $\{[\Lambda_i^L, \Lambda_i^R] : i = 1, \ldots, J\}$ be a set of $\eta$-brackets for covering $\mathcal{F}_0$ with $J \leq e^{c/\eta}$ for some constant $c$ by Theorem 2.7.5 of Van der Vaart and Wellner (1996). Then, for $\mathcal{F}_1$, we can construct a set of brackets $\{[f_i^L(X), f_i^R(X)] : i = 1, \ldots, J\}$ as follows:

$$\begin{aligned} f_i^L(X) = \sum_{j=1}^{K-1} W(T_{K,j}) \Bigg[ & \frac{\Lambda_i^L(T_{K,j}) \Delta\Lambda_0(T_{K,j+1})}{\Lambda_i^R(T_{K,j+1}) - \Lambda_i^L(T_{K,j})} \\ & - \frac{\Lambda_i^R(T_{K,j}) \Delta\Lambda_0(T_{K,j})}{\max\{\Lambda_i^L(T_{K,j}) - \Lambda_i^R(T_{K,j-1}), \Delta\Lambda_0(T_{K,j}) - 2\delta_0\}} \Bigg] \\ + W(T_{K,K}) \Bigg[ & \Lambda_i^L(T_{K,K}) \\ & - \frac{\Lambda_i^R(T_{K,K}) \Lambda_0(T_{K,K})}{\max\{\Lambda_i^L(T_{K,K}) - \Lambda_i^R(T_{K,K-1}), \Delta\Lambda_0(T_{K,K}) - 2\delta_0\}} \Bigg] \end{aligned}$$

and

$$\begin{aligned} f_i^R(X) = \sum_{j=1}^{K-1} & W(T_{K,j}) \\ & \times \Bigg[ \frac{\Lambda_i^R(T_{K,j}) \Delta\Lambda_0(T_{K,j+1})}{\max\{\Lambda_i^L(T_{K,j+1}) - \Lambda_i^R(T_{K,j}), \Delta\Lambda_0(T_{K,j+1}) - 2\delta_0\}} \\ & \qquad - \frac{\Lambda_i^L(T_{K,j}) \Delta\Lambda_0(T_{K,j})}{\Lambda_i^R(T_{K,j}) - \Lambda_i^L(T_{K,j-1})} \Bigg] \\ & + W(T_{K,K}) \Bigg[ \Lambda_i^R(T_{K,K}) - \frac{\Lambda_i^L(T_{K,K}) \Lambda_0(T_{K,K})}{\Lambda_i^R(T_{K,K}) - \Lambda_i^L(T_{K,K-1})} \Bigg]. \end{aligned}$$

It can be shown that

$$P\{f_i^R(X) - f_i^L(X)\}^2 \leq c_1 d_1^2(\Lambda_i^R, \Lambda_i^L)$$

for some constant $c_1$ and for any $\Lambda \in \mathcal{F}_0$, $P f_\Lambda^2(X) \leq c_2 d_1^2(\Lambda, \Lambda_0)$ for some constant $c_2$. Hence, $\mathcal{F}_1$ is a P-Donsker class, and it follows from Lemma 2 and Corollary 2.3.12 of Van der Vaart and Wellner (1996) that $I_{1n}^* = o_p(1)$.

Next, we show the convergence of $I_{2n}$. Set $W_0 = W \circ \Lambda_0^{-1}$. Then, from Lemma 1, we can rewrite $I_{2n}$ as

$$\begin{aligned} I_{2n} = \sqrt{n} P_n \Bigg[ &\sum_{j=1}^{K-1} \{W_0(\Lambda_0(T_{K,j})) - W_0(\hat{\Lambda}_n(T_{K,j}))\} \hat{\Lambda}_n(T_{K,j}) \\ &\times \left\{ \frac{\Delta N(T_{K,j+1})}{\Delta \hat{\Lambda}_n(T_{K,j+1})} - \frac{\Delta N(T_{K,j})}{\Delta \hat{\Lambda}_n(T_{K,j})} \right\} \\ &+ \{W_0(\Lambda_0(T_{K,K})) - W_0(\hat{\Lambda}_n(T_{K,K}))\} \hat{\Lambda}_n(T_{K,K}) \\ &\times \left\{ 1 - \frac{\Delta N(T_{K,K})}{\Delta \hat{\Lambda}_n(T_{K,K})} \right\} \Bigg] \\ = \Delta_{1n} + \Delta_{2n}, \end{aligned}$$

where

$$\begin{aligned} \Delta_{1n} = \sqrt{n}(P_n - P) \Bigg[ &\sum_{j=1}^{K-1} \{W_0(\Lambda_0(T_{K,j})) - W_0(\hat{\Lambda}_n(T_{K,j}))\} \hat{\Lambda}_n(T_{K,j}) \\ &\times \left\{ \frac{\Delta N(T_{K,j+1})}{\Delta \hat{\Lambda}_n(T_{K,j+1})} - \frac{\Delta N(T_{K,j})}{\Delta \hat{\Lambda}_n(T_{K,j})} \right\} \\ &+ \{W_0(\Lambda_0(T_{K,K})) - W_0(\hat{\Lambda}_n(T_{K,K}))\} \hat{\Lambda}_n(T_{K,K}) \\ &\times \left\{ 1 - \frac{\Delta N(T_{K,K})}{\Delta \hat{\Lambda}_n(T_{K,K})} \right\} \Bigg] \end{aligned}$$

and

$$\begin{aligned} \Delta_{2n} = \sqrt{n} P \Bigg[ &\sum_{j=1}^{K-1} \{W_0(\Lambda_0(T_{K,j})) - W_0(\hat{\Lambda}_n(T_{K,j}))\} \hat{\Lambda}_n(T_{K,j}) \\ &\times \left\{ \frac{\Delta N(T_{K,j+1})}{\Delta \hat{\Lambda}_n(T_{K,j+1})} - \frac{\Delta N(T_{K,j})}{\Delta \hat{\Lambda}_n(T_{K,j})} \right\} \\ &+ \{W_0(\Lambda_0(T_{K,K})) - W_0(\hat{\Lambda}_n(T_{K,K}))\} \hat{\Lambda}_n(T_{K,K}) \\ &\times \left\{ 1 - \frac{\Delta N(T_{K,K})}{\Delta \hat{\Lambda}_n(T_{K,K})} \right\} \Bigg]. \end{aligned}$$

Let $\Delta_{1n}^*$ and $\Delta_{2n}^*$ denote the versions of $\Delta_{1n}$ and $\Delta_{2n}$ obtained by replacing $\hat{\Lambda}_n$ with $\hat{\Lambda}_n^*$, respectively. Set

$$h_\Lambda(X) = \sum_{j=1}^{K-1} \{W_0(\Lambda_0(T_{K,j})) - W_0(\Lambda(T_{K,j}))\} \Lambda(T_{K,j})$$

$$\times \left\{ \frac{\Delta N(T_{K,j+1})}{\Delta \Lambda(T_{K,j+1})} - \frac{\Delta N(T_{K,j})}{\Delta \Lambda(T_{K,j})} \right\}$$
$$+ \{W_0(\Lambda_0(T_{K,K})) - W_0(\Lambda(T_{K,K}))\}\Lambda(T_{K,K})$$
$$\times \left\{ 1 - \frac{\Delta N(T_{K,K})}{\Delta \Lambda(T_{K,K})} \right\}$$

and

$$\mathcal{F}_2 = \{h_\Lambda(X) : \Lambda \in \mathcal{F}_0\}.$$

Note that the uniform covering entropy for $\mathcal{F}_0$ is bounded by $c/\eta$ for some constant $c$ from Theorem 2.7.5 of Van der Vaart and Wellner (1996). Since $W_0$ is a bounded Lipschitz function, it can be shown that for $\Lambda_1, \Lambda_2 \in \mathcal{F}_0$,

$$P\{(h_{\Lambda_1}(X) - h_{\Lambda_2}(X))^2\} \le c_3 d_1^2(\Lambda_1, \Lambda_2)$$

for some constant $c_3$ and for any $\Lambda \in \mathcal{F}_0$,

$$P(h_\Lambda^2(X)) \le c_4 d_1^2(\Lambda, \Lambda_0)$$

for some constant $c_4$. Hence, the uniform entropy for $\mathcal{F}_2$ is bounded by $c/\eta$, and then $\mathcal{F}_2$ is a P-Donsker class from Theorem 2.5.2 of Van der Vaart and Wellner (1996). Since $d_1(\hat{\Lambda}_n^*, \Lambda_0) \to_p 0$, it follows from the uniform asymptotic equicontinuity of the empirical process [Van der Vaart and Wellner (1996), pages 168–171] that $\Delta_{1n}^* = o_p(1)$. Then, we have $\Delta_{1n} = o_p(1)$, since $P\{\Delta_{1n} \ne \Delta_{1n}^*\} < \varepsilon$.

For $\Delta_{2n}^*$, since $W_0$ is a bounded Lipschitz function, it follows that

$$|\Delta_{2n}^*| \le c_5 \sqrt{n} d_1^2(\hat{\Lambda}_n^*, \Lambda_0),$$

where $c_5$ is a constant. This shows, from Lemma 2 and $P(\hat{\Lambda}_n \ne \hat{\Lambda}_n^*) < \varepsilon$, that $\Delta_{2n} = o_p(1)$.

For $I_{4n}$, we let $I_{4n}^*$ denote the version of $I_{4n}$ obtained by replacing $\hat{\Lambda}_n$ with $\hat{\Lambda}_n^*$, and let

$$\mathcal{F}_3 = \{g_\Lambda(X) - g_{\Lambda_0}(X) : \Lambda \in \mathcal{F}_0\}.$$

We can use the same techniques as those used for proving the convergence of $I_{1n}$ to show that $\mathcal{F}_3$ is P-Donsker and $P\{g_\Lambda(X) - g_{\Lambda_0}(X)\}^2 \le c_6 d_1^2(\Lambda, \Lambda_0)$ for some constant $c_6$, and hence $I_{4n}^* = o_p(1)$, which completes the proof of the theorem.

6.2. *Proof of Theorem 3.1.* (i) To obtain the asymptotic distribution of $\mathbf{U}_n$, we first note that $U_n^{(l)}$ can rewritten as

$$U_n^{(l)} = U_{1n}^{(l)} - \sqrt{\frac{n}{n_l}} U_{2n}^{(l)} + U_{3n}^{(l)} + U_{4n}^{(l)} + U_{5n}^{(l)} + U_{6n}^{(l)}, \tag{6.15}$$

where, for $l = 1, \ldots, k$,

$$U_{1n}^{(l)} = \sqrt{n} P\Bigg[\sum_{j=1}^{K-1} W(T_{K,j})\hat{\Lambda}_n(T_{K,j})\bigg\{\frac{\Delta\Lambda_0(T_{K,j+1})}{\Delta\hat{\Lambda}_n(T_{K,j+1})} - \frac{\Delta\Lambda_0(T_{K,j})}{\Delta\hat{\Lambda}_n(T_{K,j})}\bigg\} + W(T_{K,K})\hat{\Lambda}_n(T_{K,K})\bigg\{1 - \frac{\Delta\Lambda_0(T_{K,K})}{\Delta\hat{\Lambda}_n(T_{K,K})}\bigg\}\Bigg],$$

$$U_{2n}^{(l)} = \sqrt{n_l} P\Bigg[\sum_{j=1}^{K-1} W(T_{K,j})\hat{\Lambda}_{n_l}(T_{K,j})\bigg\{\frac{\Delta\Lambda_0(T_{K,j+1})}{\Delta\hat{\Lambda}_{n_l}(T_{K,j+1})} - \frac{\Delta\Lambda_0(T_{K,j})}{\Delta\hat{\Lambda}_{n_l}(T_{K,j})}\bigg\} + W(T_{K,K})\hat{\Lambda}_{n_l}(T_{K,K})\bigg\{1 - \frac{\Delta\Lambda_0(T_{K,K})}{\Delta\hat{\Lambda}_{n_l}(T_{K,K})}\bigg\}\Bigg],$$

$$U_{3n}^{(l)} = \sqrt{n}(P_n - P)\Bigg[\sum_{j=1}^{K-1} W_n^{(l)}(T_{K,j})\hat{\Lambda}_n(T_{K,j}) \times \bigg\{\frac{\Delta\hat{\Lambda}_{n_l}(T_{K,j+1})}{\Delta\hat{\Lambda}_n(T_{K,j+1})} - \frac{\Delta\hat{\Lambda}_{n_l}(T_{K,j})}{\Delta\hat{\Lambda}_n(T_{K,j})}\bigg\} + W(T_{K,K})\hat{\Lambda}_n(T_{K,K}) \times \bigg\{1 - \frac{\Delta\hat{\Lambda}_{n_l}(T_{K,K})}{\Delta\hat{\Lambda}_n(T_{K,K})}\bigg\}\Bigg],$$

$$U_{4n}^{(l)} = \sqrt{n} P\Bigg[\sum_{j=1}^{K-1} \{W_n^{(l)}(T_{K,j}) - W(T_{K,j})\}\hat{\Lambda}_n(T_{K,j}) \times \bigg\{\frac{\Delta\hat{\Lambda}_{n_l}(T_{K,j+1})}{\Delta\hat{\Lambda}_n(T_{K,j+1})} - \frac{\Delta\hat{\Lambda}_{n_l}(T_{K,j})}{\Delta\hat{\Lambda}_n(T_{K,j})}\bigg\} + \{W_n^{(l)}(T_{K,K}) - W(T_{K,K})\}\hat{\Lambda}_n(T_{K,K}) \times \bigg\{1 - \frac{\Delta\hat{\Lambda}_{n_l}(T_{K,K})}{\Delta\hat{\Lambda}_n(T_{K,K})}\bigg\}\Bigg],$$

$$U_{5n}^{(l)} = \sqrt{n} P\Bigg[\sum_{j=1}^{K-1} W(T_{K,j})\{\hat{\Lambda}_n(T_{K,j}) - \hat{\Lambda}_{n_l}(T_{K,j})\} \times \bigg\{\frac{\Delta\hat{\Lambda}_{n_l}(T_{K,j+1}) - \Delta\Lambda_0(T_{K,j+1})}{\Delta\hat{\Lambda}_n(T_{K,j+1})} - \frac{\Delta\hat{\Lambda}_{n_l}(T_{K,j}) - \Delta\Lambda_0(T_{K,j})}{\Delta\hat{\Lambda}_n(T_{K,j})}\bigg\}$$

$$+ W(T_{K,K})\hat{\Lambda}_n(T_{K,K})\left\{-\frac{\Delta\hat{\Lambda}_{n_l}(T_{K,K}) - \Delta\Lambda_0(T_{K,K})}{\Delta\hat{\Lambda}_n(T_{K,K})}\right\}\Bigg]$$

and

$$\begin{aligned}
U_{6n}^{(l)} = \sqrt{n}P\Bigg[\sum_{j=1}^{K-1} W(T_{K,j})\hat{\Lambda}_{n_l}(T_{K,j})\bigg\{&(\Delta\hat{\Lambda}_{n_l}(T_{K,j+1}) - \Delta\hat{\Lambda}_0(T_{K,j+1})) \\
&\times\left(\frac{1}{\Delta\hat{\Lambda}_n(T_{K,j+1})} - \frac{1}{\Delta\hat{\Lambda}_{n_l}(T_{K,j+1})}\right) \\
&-(\Delta\hat{\Lambda}_{n_l}(T_{K,j}) - \Delta\hat{\Lambda}_0(T_{K,j})) \\
&\times\left(\frac{1}{\Delta\hat{\Lambda}_n(T_{K,j})} - \frac{1}{\Delta\hat{\Lambda}_{n_l}(T_{K,j})}\right)\bigg\} \\
- W(T_{K,K})\hat{\Lambda}_{n_l}(T_{K,K})\{\Delta\hat{\Lambda}_{n_l}(T_{K,K}) - \Delta\Lambda_0(T_{K,K})\}& \\
\times\left\{\frac{1}{\Delta\hat{\Lambda}_n(T_{K,K})} - \frac{1}{\Delta\hat{\Lambda}_{n_l}(T_{K,K})}\right\}\Bigg].&
\end{aligned}$$

From the proof of Theorems 2.1, we have, for $l = 1, \ldots, k$,

$$U_{1n}^{(l)} = Y_n + o_p(1)$$

and

$$U_{2n}^{(l)} = Y_n^{(l)} + o_p(1),$$

where

$$\begin{aligned}
Y_n = \sqrt{n}(P_n - P)\Bigg[&\sum_{j=1}^{K-1} W(T_{K,j})\Lambda_0(T_{K,j})\left\{\frac{\Delta N(T_{K,j+1})}{\Delta\Lambda_0(T_{K,j+1})} - \frac{\Delta N(T_{K,j})}{\Delta\Lambda_0(T_{K,j})}\right\} \\
&+ W(T_{K,K})\Lambda_0(T_{K,K})\left\{1 - \frac{\Delta N(T_{K,K})}{\Delta\Lambda_0(T_{K,K})}\right\}\Bigg]
\end{aligned}$$

and

$$\begin{aligned}
Y_n^{(l)} = \sqrt{n_l}(P_{n_l} - P)\Bigg[&\sum_{j=1}^{K-1} W(T_{K,j})\Lambda_0(T_{K,j})\left\{\frac{\Delta N(T_{K,j+1})}{\Delta\Lambda_0(T_{K,j+1})} - \frac{\Delta N(T_{K,j})}{\Delta\Lambda_0(T_{K,j})}\right\} \\
&+ W(T_{K,K})\Lambda_0(T_{K,K})\left\{1 - \frac{\Delta N(T_{K,K})}{\Delta\Lambda_0(T_{K,K})}\right\}\Bigg],
\end{aligned}$$

where $P_{n_l}f = \frac{1}{n_l}\sum_{i\in S_l} f(X_i)$ and $S_l$ denotes the set of indices for subjects in group $l$, $l = 1, \ldots, k$. Evidently, $Y_n^{(l)}$'s are independent and identically distributed, and $\sqrt{n}Y_n = \sum_{l=1}^{k} \sqrt{n_l} Y_n^{(l)}$. Set $Z_n^{(l)} = Y_n - \sqrt{\frac{n}{n_l}} Y_n^{(l)}$, $l = 1, \ldots, k$ and $\mathbf{Z}_n = (Z_n^{(1)}, \ldots, Z_n^{(k)})^T$. Then,

$$Z_n^{(l)} = \sum_{i=1}^{k} \sqrt{\frac{n_i}{n}} Y_n^{(i)} - \sqrt{\frac{n}{n_l}} Y_n^{(l)}, \qquad l = 1, \ldots, k,$$

and so

$$\mathbf{Z}_n = \mathbf{\Gamma}_n \mathbf{Y}_n = \mathbf{\Gamma} \mathbf{Y}_n + o_p(1),$$

where

$$\mathbf{\Gamma} = \begin{pmatrix} \sqrt{p_1} - \dfrac{1}{\sqrt{p_1}} & \sqrt{p_2} & \cdots & \sqrt{p_k} \\ \sqrt{p_1} & \sqrt{p_2} - \dfrac{1}{\sqrt{p_2}} & \cdots & \sqrt{p_k} \\ \cdots & \cdots & \cdots & \cdots \\ \sqrt{p_1} & \sqrt{p_2} & \cdots & \sqrt{p_k} - \dfrac{1}{\sqrt{p_k}} \end{pmatrix}$$

and

$$\mathbf{Y}_n = (Y_n^{(1)}, \ldots, Y_n^{(k)})^T$$

converges in distribution to $\mathbf{Y}_w$, having a $k$-dimensional normal distribution with mean vector $\mathbf{0}$ and covariance matrix $\operatorname{diag}(\sigma_1^2, \ldots, \sigma_k^2)$, where $\sigma_l^2$'s are given in the statement of the theorem. Thus, we have $\mathbf{Z}_n$ converging in distribution to a random variable $\mathbf{U}_w$ that has a normal distribution $N(\mathbf{0}, \mathbf{\Sigma}_{\mathbf{U}_w})$, where $\mathbf{\Sigma}_{\mathbf{U}_w}$ is given in (3.4) of Theorem 3.1.

Now, we need to show that $U_{3n}^{(l)}$, $U_{4n}^{(l)}$, $U_{5n}^{(l)}$ and $U_{6n}^{(l)}$ all converge in probability to 0, $l = 1, \ldots, k$. Let $U_{3n}^{(l)*}$, $U_{4n}^{(l)*}$, $U_{5n}^{(l)*}$ and $U_{6n}^{(l)*}$ denote the version of $U_{3n}^{(l)}$, $U_{4n}^{(l)}$, $U_{5n}^{(l)}$ and $U_{6n}^{(l)}$ obtained by replacing $\hat{\Lambda}_n$ with $\hat{\Lambda}_n^*$ and $\hat{\Lambda}_{n_l}$ with $\hat{\Lambda}_{n_l}^*$, respectively. Then, to prove that $U_{3n}^{(l)}$, $U_{4n}^{(l)}$, $U_{5n}^{(l)}$ and $U_{6n}^{(l)}$ all converge in probability to 0, $l = 1, \ldots, k$, it is sufficient to show that $U_{3n}^{(l)*}$, $U_{4n}^{(l)*}$, $U_{5n}^{(l)*}$, and $U_{6n}^{(l)*}$ all converge in probability to 0, $l = 1, \ldots, k$.

For $U_{3n}^{(l)*}$, set

$$\mathcal{G} = \{\xi : [0, \tau] \to [0, b]\},$$

where $b$ is the uniform upper bound of weight process $W_n^{(l)}$ $(l = 1, \ldots, k)$,

$$\psi_{\Lambda_1, \Lambda_2}(\xi, X) = \sum_{j=1}^{K-1} \xi(T_{K,j}) \Lambda_1(T_{K,j}) \left\{ \frac{\Delta\Lambda_2(T_{K,j+1})}{\Delta\Lambda_1(T_{K,j+1})} - \frac{\Delta\Lambda_2(T_{K,j})}{\Delta\Lambda_1(T_{K,j})} \right\} + \xi(T_{K,K}) \Lambda_1(T_{K,K}) \left\{ 1 - \frac{\Delta\Lambda_2(T_{K,K})}{\Delta\Lambda_1(T_{K,K})} \right\},$$

and, for $\xi \in \mathcal{G}$,

$$\Psi_\delta(\xi) = \{\psi_{\Lambda_1,\Lambda_2}(\xi, X) : \Lambda_1, \Lambda_2 \in \mathcal{F}_\delta\}.$$

Note that it follows from Theorem 2.7.5 of Van der Vaart and Wellner (1996) that

$$N_{[\cdot]}(\eta, \mathcal{F}_\delta, L_2(P)) \le e^{c_1/\eta}$$

for some constant $c_1$. Then, we have

$$N_{[\cdot]}(\eta, \Psi_\delta(\xi), L_2(P)) \le e^{2c_1/\eta}.$$

It can be easily shown that $|\psi_{\Lambda_1,\Lambda_2}(\xi, X)| \le \psi(X)$, where $\psi(X) = c_2 K$, and $P\psi^2_{\Lambda_1,\Lambda_2}(\xi, X) \le c_3\delta^2$, where $c_2$ and $c_3$ are universal constants for $\xi$. Thus,

$$J_{[\cdot]}(\delta, \Psi_\delta(\xi), L_2(P)) = \int_0^\delta \sqrt{1 + \log N_{[\cdot]}(\eta\|\psi\|_{P,2}, \Psi_\delta(\xi), L_2(P))}\, d\eta \le c_4\delta^{1/2}$$

for some constant universal $c_4$. Hence, from Theorem 2.14.2 of Van der Vaart and Wellner (1996), we have

$$\begin{aligned} &E^*\left\{\sup_{\psi_{\Lambda_1,\Lambda_2}(\xi,X)\in\Psi_\delta(\xi)} |\sqrt{n}(P_n - P)\psi_{\Lambda_1,\Lambda_2}(\xi, X)|\right\} \\ &\qquad \le c_5[J_{[\cdot]}(\delta, \Psi_\delta(\xi), L_2(P)) + \sqrt{n}P\psi\{\psi > \sqrt{n}a(\delta)\}], \end{aligned}$$

where $c_5$ is a universal constant and

$$a(\delta) = \delta\|\psi\|_{P,2} / \sqrt{1 + \log N_{[\cdot]}(\delta\|\psi\|_{P,2}, \Psi_\delta(\xi), L_2(P))}.$$

Then, it can be easily shown that

$$\limsup_{n\to\infty} E^*\left\{\sup_{\psi_{\Lambda_1,\Lambda_2}(\xi,X)\in\Psi_\delta(\xi)} |\sqrt{n}(P_n - P)\psi_{\Lambda_1,\Lambda_2}(\xi, X)|\right\} \le c_6\delta^{1/2}$$

for some universal constant $c_6$. It follows from $d_1(\hat\Lambda_n, \hat\Lambda_{n_l}) \xrightarrow{a.s.} 0$ that

$$\limsup_{n\to\infty} E|\sqrt{n}(P_n - P)\psi_{\hat\Lambda^*_n,\hat\Lambda^*_{n_l}}(W_n^{(l)}, X)| \le c_6\delta^{1/2}.$$

Letting $\delta \to 0$, we have

$$\lim_{n\to\infty} E|\sqrt{n}(P_n - P)\psi_{\hat\Lambda^*_n,\hat\Lambda^*_{n_l}}(W_n^{(l)}, X)| = 0,$$

which yields $U_{3n}^{(l)*} = o_p(1)$.

For $U_{4n}^{(l)*}$, we note that

$$\begin{aligned}|U_{4n}^{(l)*}| &\le c_7\Bigg[\sqrt{n}P\Bigg\{\sum_{j=1}^{K}|W_n^{(l)}(T_{K,j-1}) - W(T_{K,j-1})| \\ &\qquad\qquad \times |\hat{\Lambda}_n^*(T_{K,j}) - \hat{\Lambda}_{n_l}^*(T_{K,j})|\Bigg\} \\ &\quad + \sqrt{n}P\Bigg\{\sum_{j=1}^{K}|W_n^{(l)}(T_{K,j}) - W(T_{K,j})||\hat{\Lambda}_n^*(T_{K,j}) - \hat{\Lambda}_{n_l}^*(T_{K,j})|\Bigg\} \\ &\quad + \sqrt{n}P\Bigg\{\sum_{j=1}^{K}|W_n^{(l)}(T_{K,j}) - W(T_{K,j})| \\ &\qquad\qquad \times |\hat{\Lambda}_n^*(T_{K,j-1}) - \hat{\Lambda}_{n_l}^*(T_{K,j-1})|\Bigg\}\Bigg] \\ &= c_7(A_{1n}^{(l)} + A_{2n}^{(l)} + A_{3n}^{(l)})\end{aligned}$$

for some constant $c_7$, where

$$\begin{aligned}A_{1n}^{(l)} &= \sqrt{n}\iint |W_n^{(l)}(u) - W(u)||\hat{\Lambda}_n^*(v) - \hat{\Lambda}_{n_l}^*(v)|\, d\mu_2(u,v) \\ &\le \sqrt{n}\iint |W_n^{(l)}(u) - W(u)||\hat{\Lambda}_n^*(v) - \Lambda_0(v)|\, d\mu_2(u,v) \\ &\quad + \sqrt{n}\iint |W_n^{(l)}(u) - W(u)||\hat{\Lambda}_{n_l}^*(v) - \Lambda_0(v)|\, d\mu_2(u,v), \\ A_{2n}^{(l)} &= \sqrt{n}\int_0^\tau |W_n^{(l)}(t) - W(t)||\hat{\Lambda}_n^*(t) - \hat{\Lambda}_{n_l}^*(t)|\, d\mu_1(t) \\ &\le \sqrt{n}\int_0^\tau |W_n^{(l)}(t) - W(t)||\hat{\Lambda}_n^*(t) - \Lambda_0(t)|\, d\mu_1(t) \\ &\quad + \sqrt{n}\int_0^\tau |W_n^{(l)}(t) - W(t)||\hat{\Lambda}_{n_l}^*(t) - \Lambda_0(t)|\, d\mu_1(t)\end{aligned}$$

and

$$\begin{aligned}A_{3n}^{(l)} &= \sqrt{n}\iint |W_n^{(l)}(v) - W(v)||\hat{\Lambda}_n^*(u) - \hat{\Lambda}_{n_l}^*(u)|\, d\mu_2(u,v) \\ &\le \sqrt{n}\iint |W_n^{(l)}(v) - W(v)||\hat{\Lambda}_n^*(u) - \Lambda_0(u)|\, d\mu_2(u,v) \\ &\quad + \sqrt{n}\iint |W_n^{(l)}(v) - W(v)||\hat{\Lambda}_{n_l}^*(u) - \Lambda_0(u)|\, d\mu_2(u,v).\end{aligned}$$

Using the Cauchy–Schwarz inequality, we have

$$\begin{aligned}
\sqrt{n}\iint &|W_n^{(l)}(u) - W(u)||\hat{\Lambda}_n^*(v) - \Lambda_0(v)|\, d\mu_2(u,v) \\
&\le c_8\sqrt{n}\left\{\int_0^\tau (W_n^{(l)}(t) - W(t))^2\, d\mu_1(t)\right\}^{1/2} \\
&\quad \times \left\{\int_0^\tau (\hat{\Lambda}_n^*(t) - \Lambda_0(t))^2\, d\mu_1(t)\right\}^{1/2} \\
&\longrightarrow 0,
\end{aligned}$$

in probability, where $c_8$ is a constant, since

$$\left[\int_0^\tau \{\hat{\Lambda}_n^*(t) - \Lambda_0(t)\}^2\, d\mu_1(t)\right]^{1/2} = O_p(n^{-1/3}).$$

Similarly, we have

$$\sqrt{n}\iint |W_n^{(l)}(u) - W(u)||\hat{\Lambda}_{n_l}^*(v) - \Lambda_0(v)|\, d\mu_2(u,v) = o_p(1).$$

Thus, $A_{1n}^{(l)} = o_p(1)$. Similarly, we have $A_{2n}^{(l)} = o_p(1)$ and $A_{3n}^{(l)} = o_p(1)$. Hence, $U_{4n}^{(l)*} = o_p(1)$, $l = 1, \ldots, k$.

For $U_{5n}^{(l)*}$ and $U_{6n}^{(l)*}$, we note that

$$\begin{aligned}
|U_{5n}^{(l)*}| &\le c_9\left[\sqrt{n}P\left\{\sum_{j=1}^K |\hat{\Lambda}_n^*(T_{K,j}) - \Lambda_0(T_{K,j})|^2\right\}\right. \\
&\qquad \left. + \sqrt{n}P\left\{\sum_{j=1}^K |\hat{\Lambda}_{n_l}^*(T_{K,j}) - \Lambda_0(T_{K,j})|^2\right\}\right] \\
&= c_9\{\sqrt{n}d_1^2(\hat{\Lambda}_n^*, \Lambda_0) + \sqrt{n}d_1^2(\hat{\Lambda}_{n_l}^*, \Lambda_0)\}
\end{aligned}$$

and

$$|U_{6n}^{(l)*}| \le c_{10}\{\sqrt{n}d_1^2(\hat{\Lambda}_n^*, \Lambda_0) + \sqrt{n}d_1^2(\hat{\Lambda}_{n_l}^*, \Lambda_0)\}$$

for some constants $c_9$ and $c_{10}$. Hence, $U_{5n}^{(l)*} = o_p(1)$ and $U_{6n}^{(l)*} = o_p(1)$, $l = 1, \ldots, k$. Therefore, the proof of part (i) is complete.

(ii) We note that $V_n^{(l)} = U_n^{(1,l)} - U_n^{(l)}$, $l = 2, \ldots, k$, where $U_n^{(1,l)}$ is defined as $U_n^{(1)}$ by replacing $W_n^{(1)}$ with $W_n^{(l)}$ for $l = 2, \ldots, k$. Then, it follows from (i) that

$$V_n^{(l)} = -\sqrt{\frac{n}{n_1}}Y_n^{(1)} + \sqrt{\frac{n}{n_l}}Y_n^{(l)} + o_p(1)$$

for $l = 2, \ldots, k$, and so

$$\mathbf{V}_n = \mathbf{H}_n \mathbf{Y}_n + o_p(1) = \mathbf{H}\mathbf{Y}_n + o_p(1),$$

where $\mathbf{H}_n$ and $\mathbf{H}$ are given in the theorem. This completes the proof of part (ii).

(iii) To show that $\hat{\sigma}_l^2 - \sigma_w^2 = o_p(1)$ for $l = 1, \ldots, k$, we set

$$\begin{aligned}\phi(\xi, \Lambda, X) &= \sum_{j=1}^{K-1} \xi(T_{K,j})\Lambda(T_{K,j})\left\{\frac{\Delta N(T_{K,j+1})}{\Delta\Lambda(T_{K,j+1})} - \frac{\Delta N(T_{K,j})}{\Delta\Lambda(T_{K,j})}\right\} \\ &\quad + \xi(T_{K,K})\Lambda(T_{K,K})\left\{1 - \frac{\Delta N(T_{K,K})}{\Delta\Lambda(T_{K,K})}\right\}.\end{aligned}$$

Then, $\sigma_w^2 = P\phi^2(W, \Lambda_0, X)$ and $\hat{\sigma}_l^2 = P_n\phi^2(W_n^{(l)}, \hat{\Lambda}_n, X)$. Note that

$$\begin{aligned}\hat{\sigma}_l^2 - \sigma_w^2 &= P_n\{\phi^2(W_n^{(l)}, \hat{\Lambda}_n, X) - \phi^2(W_n^{(l)}, \Lambda_0, X)\} \\ &\quad + P_n\{\phi^2(W_n^{(l)}, \Lambda_0, X) - \phi^2(W, \Lambda_0, X)\} \\ &\quad + (P_n - P)\phi^2(W, \Lambda_0, X).\end{aligned}$$

It can be easily shown that

$$P_n\{\phi^2(W_n^{(l)}, \hat{\Lambda}_n, X) - \phi^2(W_n^{(l)}, \Lambda_0, X)\} = o_p(1)$$

and

$$(P_n - P)\phi^2(W_0, \Lambda_0, X) = o_p(1).$$

Since it follows from Condition C that

$$\begin{aligned}&|\phi(W_n^{(l)}, \Lambda_0, X) - \phi(W, \Lambda_0, X)| \\ &\quad = |\phi(W_n^{(l)} - W, \Lambda_0, X)| \\ &\quad \le b_1\{1 + N(T_{K,K})\}\sum_{j=1}^{K}|W_n^{(l)}(T_{K,j}) - W(T_{K,j})|\end{aligned}$$

with probability 1 for some constant $b_1$ and

$$\begin{aligned}|\phi(W_n^{(l)}, \Lambda_0, X) + \phi(W, \Lambda_0, X)| &= |\phi(W_n^{(l)} + W, \Lambda_0, X)| \\ &\le b_2 K\{1 + N(T_{K,K})\}\end{aligned}$$

with probability 1 for some constant $b_2$, then we have, from the Cauchy–Schwarz inequality, Conditions $P\{K \le K_0\} = 1$ and D, and (3.6)

$$\begin{aligned}&E|\phi^2(W_n^{(l)}, \Lambda_0, X_i) - \phi^2(W, \Lambda_0, X_i)| \\ &\quad \le b_3 E\left[\{1 + N_i(T_{K_i,K_i})\}^2\left\{\sum_{j=1}^{K_i}|W_n^{(l)}(T_{K_i,j}) - W(T_{K_i,j})|\right\}\right]\end{aligned}$$

$$
\begin{aligned}
&\le b_3[E\{1+N_i(T_{K_i,K_i})\}^4]^{1/2}\left[E\left\{\sum_{j=1}^{K_i}|W_n^{(l)}(T_{K_i,j})-W(T_{K_i,j})|\right\}^2\right]^{1/2}\\
&\le b_4\max_{1\le i\le n}\left[E\left\{\sum_{j=1}^{K_i}|W_n^{(l)}(T_{K_i,j})-W(T_{K_i,j})|^2\right\}\right]^{1/2}\\
&\longrightarrow 0,
\end{aligned}
$$

where $b_3$ and $b_4$ are finite positive constants, which completes the proof of part (iii).

REMARK 3. The monotonicity assumption of the weight process required by Zhang (2006) can be removed by using the same techniques as those used here.

**Acknowledgments.** The authors are very grateful to the editor, Professor Susan Murphy and two referees for their helpful comments and suggestions that greatly improved the paper. They thank Ying Zhang and Minggen Lu for supplying the R-code for the modified iterative convex minorant algorithm presented by Wellner and Zhang (2000).

DEPARTMENT OF MATHEMATICS AND STATISTICS
MCMASTER UNIVERSITY
1280 MAIN STREET WEST
HAMILTON, ONTARIO
CANADA L8S 4K1
E-MAIL: bala@mcmaster.ca

DEPARTMENT OF APPLIED MATHEMATICS
HONG KONG POLYTECHNIC UNIVERSITY
6/F, CORE J, STANLEY HO BUILDING
HUNG HOM, KOWLOON
HONG KONG
E-MAIL: xingqiu.zhao@polyu.edu.hk